\documentclass[12pt]{amsart}

\newtheorem{theorem}{Theorem}[section]
\newtheorem{lemma}[theorem]{Lemma}
\newtheorem{corollary}[theorem]{Corollary}
\newtheorem{proposition}[theorem]{Proposition}

\newtheorem{conjecture}[theorem]{Conjecture}

\theoremstyle{plain}
\newtheorem{definition}[theorem]{Definition}

\newtheorem{remark}[theorem]{Remark}

\newtheorem{question}[theorem]{Question}
\newtheorem{setting}[theorem]{Setting}

\theoremstyle{definition}

\theoremstyle{remark}

\numberwithin{equation}{section}

\usepackage{fancyhdr}
\usepackage{amscd}
\usepackage{amsmath}
\usepackage{amsthm}
\usepackage{amssymb}

\usepackage{mathrsfs}

\begin{document}

\title[Chern character and semi-regularity map]{Chern character, semi-regularity map and obstructions}

\author{Sen Yang}

\address{Shing-Tung Yau Center of Southeast University \\ 
Southeast University \\
Nanjing, China\\
}
\address{School of Mathematics \\  Southeast University \\
Nanjing, China\\
}

\email{101012424@seu.edu.cn}

\date{}

\begin{abstract}
Using Chern character, we construct a natural transformation from the local Hilbert functor to a functor of Artin rings defined from Hochschild homology, which enables us to reconstruct the semi-regularity map and the infinitesimal Abel-Jacobi map.

Combining that construction of the semi-regularity map with obstruction theory of functors of Artin rings, we give a different proof of a theorem of Bloch stating that the semi-regularity map annihilates certain obstructions to embedded deformations of a closed subvariety which is a locally complete intersection.
\end{abstract}

\maketitle


 \section{Introduction}
Let $X$ be a smooth projective variety over a field $k$ of characteristic zero and let $Y \subset X$ be a closed subvariety of codimension $p$, which is a locally complete intersection. In the seminal work \cite{Bloch1}, Bloch used duality theory to construct the semi-regularity map 
\begin{equation*}
   H^{1}(N_{Y/X}) \to H^{p+1}(\Omega_{X/k}^{p-1}),
 \end{equation*}
which generalized the semi-regularity map introduced by Severi \cite{Severi} (for curves on a surface) and Kodaira-Spencer \cite{KS} (for divisors on a projective variety).

This paper is motivated by the semi-regularity conjecture by Bloch:
\begin{conjecture} \label{Conjecture:main}
The semi-regularity map $H^{1}(N_{Y/X}) \to H^{p+1}(\Omega_{X/k}^{p-1})$ annihilates every obstructions to embedded deformations of $Y$ in $X$, where $N_{Y/X}$ is the normal bundle.
\end{conjecture}

Because of its importance in deformation theory, the semi-regularity conjecture had been intensively studied. Bloch \cite{Bloch1} made significant progress on it and proved Conjecture \ref{Conjecture:main} for certain obstructions. Since we eventually look at obstructions arising from principal small extensions (see Remark \ref{Remark: obs by principal}), we recall Bloch's theorem as follows:
\begin{theorem} [\cite{Bloch1}] \label{t:Bloch-semi-conj}
For principal small extensions
\begin{equation}
e: 0 \to (\eta) \to B \to A \to 0
\end{equation}
such that the differential $(\eta) \to \Omega_{B/k} \otimes_{B}A$ is injective (this includes curvilinear extensions), the conjecture \ref{Conjecture:main} is true.
\end{theorem}

The relation between the semi-regularity map and obstructions to deformations of subvarieties had been further studied by Clemens \cite{Clemens}, Manetti \cite{Manetti-LieOb}, Ran \cite{Ran93, Ran99} and others. Using cotangent complexes, Buchweitz and Flenner \cite{BF} generalized Bloch's construction of the semi-regularity map to perfect complexes on arbitrary complex spaces.

In the hope of relating both the infinitesimal Abel-Jacobi map 
\[
H^{0}(N_{Y/X}) \to H^{p}(X, \Omega^{p-1}_{X/k})
\]
and the semi-regularity map with a natural transformation of two deformation functors, one is interested in the following question (see page 138 of \cite{BF}):
\begin{question}  \label{Question:1}
Is it possible to find two deformation functors $\mathbb{F}_{1}$ and $\mathbb{F}_{2}$, and then construct a natural transformation of them such that

$\mathrm{(I)}$ the infinitesimal Abel-Jacobi map maps the tangent space of $\mathbb{F}_{1}$ to the tangent space of  $\mathbb{F}_{2}$;

$\mathrm{(II)}$ the semi-regularity map maps the obstruction space of $\mathbb{F}_{1}$ to the obstruction space of  $\mathbb{F}_{2}$?
\end{question}

This question provides an approach to Conjecture \ref{Conjecture:main}.  It is obvious that a natural transformation from one deformation functor to a smooth deformation functor induces a map of obstruction spaces, which annihilates every obstructions. This is the idea used by Iacono and Manetti \cite{IM}, in which they proved Conjecture \ref{Conjecture:main} under an additional assumption on the existence of extendable normal bundle. Two deformation functors used by Iacono and Manetti \cite{IM} are constructed from differential graded Lie algebras and the natural transformation of them is induced by a morphism of differential graded Lie algebras. The guiding principle of this construction is that a deformation problem over a field of characteristic zero is governed by a differential graded Lie algebra, which is known to Deligne, Drinfeld, Hitchin, Kontsevich, Lurie, Manetti, Pridham and others. It is noted that an interesting contribution to Conjecture \ref{Conjecture:main} was made by Pridham \cite{P}.

Deformation functors are functors of Artin rings satisfying certain Schlessinger type conditions, see Definition \ref{Definition: Deformation functor} in appendix. Instead of considering deformation functors, we consider functors of Artin rings and modify Question \ref{Question:1} as following:
\begin{question} \label{Question:2}
Is it possible to find two functors of Artin rings $\mathbb{F}_{1}$ and $\mathbb{F}_{2}$, and then construct a natural transformation of them such that

$\mathrm{(I)}$ the infinitesimal Abel-Jacobi map maps the tangent space of $\mathbb{F}_{1}$ to the tangent space of  $\mathbb{F}_{2}$;

$\mathrm{(II)}$ the semi-regularity map maps the obstruction space of $\mathbb{F}_{1}$ to the obstruction space of  $\mathbb{F}_{2}$?

\end{question}

While the functor $\mathbb{F}_{1}$ is expected to be the local Hilbert functor, the search of the functor $\mathbb{F}_{2}$ is open. 

Guided by Question \ref{Question:2}, this paper is organized as follows. Obstruction theory of functors of Artin rings is recalled in section 2. In the third section,  we use Hochschild homology to introduce a functor of Artin rings $\mathbb{HH}^{(p)}$ in Definition \ref{Definition: Functor-HH-final}. Moreover, we construct a natural transformation $\mathbb{T}$ from the local Hilbert functor to the functor $\mathbb{HH}^{(p)}$ in Lemma \ref{Lemma: transf-Hilb-HH-final}. Explicit computations of Hochschild homology have also been done. 

In section 4, we use the natural transformation $\mathbb{T}$ to reconstruct the infinitesimal Abel-Jacobi map and the semi-regularity map. Moreover, we answer affirmatively (after slight modification) the part (I) of Question \ref{Question:2} in Theorem \ref{Theorem: Inf'l AJ-tang-space}, and give a different proof of Theorem \ref{t:Bloch-semi-conj} in Theorem \ref{Theorem: Semi-reg obs space}. In appendix, we explain that the functor of Artin rings $\mathbb{HH}^{(p)}$ is not a deformation functor.

The functor $\mathbb{HH}^{(p)}$ in Definition \ref{Definition: Functor-HH-final} and the natural transformation $\mathbb{T}$ in Lemma \ref{Lemma: transf-Hilb-HH-final} play a key role in this paper. We record that the definitions of $\mathbb{HH}^{(p)}$ and $\mathbb{T}$ are inspired by works of Green-Griffiths \cite{GGtangentspace} and Weibel \cite{W3} as follows.

(1). Green and Griffiths asked the following question in section 7.2 of \cite{GGtangentspace} (see also Question 1.2 in \cite{Y-3}):
\begin{question}  \label{question: comparetangent}
Let $X$ be a smooth projective variety over a field $k$ of characteristic zero and let $Y \subset X$ be a closed subvariety of codimension $p$, is it possible to define a map from the tangent space $\mathrm{T_{Y}Hilb}^{p}(X)$ of the Hilbert scheme at the point $Y$ to the tangent space of the cycle group $TZ^{p}(X)$
 \[
  \mathrm{T_{Y}Hilb}^{p}(X)  \to TZ^{p}(X) ?
 \]
\end{question}
For $p=\mathrm{dim}(X)$, this question has been answered by Green-Griffiths \cite{GGtangentspace}; for general $p$, it has been answered by the author \cite{Y-3}. This question guides us to consider Chern character from Grothendieck group to negative cyclic homology.

(2). According to Weibel (see \textbf{Chern characters 3.8} of \cite{W3}), the composition of the Chern character from Grothendieck group to negative cyclic homology $K_{0}(X) \to HN_{0}(X)$ with the natural map from negative cyclic homology to Hochschild homology $HN_{0}(X) \to HH_{0}(X)$ yields Dennis trace map, which completely determines the classical Chern character of $X$ (see Proposition 4.5 of \cite{W3}). This highlights the importance of Hochschild homology and guides us to consider the natural map from Grothendieck group to 
Hochschild homology, which is the main ingredient of constructing the natural transformation $\mathbb{T}$ in Lemma \ref{Lemma: transf-Hilb-HH-final}.

\textbf{Notation}:

(1). K-theory used in this paper is Thomason-Trobaugh K-theory. 

(2). Hochschild homology and negative cyclic homology are defined over a field $k$ of characteristic zero, if not stated otherwise. 

(3). The ring of dual numbers is denoted by $k[\varepsilon]$.

(4). For $F$ a quasi-coherent sheaf, $\mathcal{H}^{*}_{Y}(F)$ denotes local cohomology sheaf, which is obtained by localizing local cohomology $H^{*}_{Y}(F)$.

\section{Obstruction theory of Functors of Artin rings}
Fucntors of Artin rings introduced by Schlessinger \cite{Schless} is an important tool in deformation theory. Fantechi-Manetti \cite{FM-Obs} developed relative obstruction theory for a morphism of functors of Artin ring $F \to G$. When $G=*$, this gives obstruction theory of a functor of Artin rings $F$ which is recalled below.

To fix notations, $k$ is a field of characteristic zero. Let $Art_{k}$ denote the category of local artinian $k$-algebra with residue field $k$, the morphisms are local homomorphisms. Let $Set_{*}$ denote the category of pointed sets. For an element $V \in Set_{*}$, $*$ denotes the chosen point of $V$; when $V$ is a vector space, the chosen point is zero.

\begin{definition}
A functor of Artin rings is a covariant functor $F: Art_{k} \to Set_{*}$ such that $F(k)=*$.
\end{definition}

\begin{definition} [Definition 1.0 of \cite{FM-Obs}] \label{Definition: Small Extension}
A small extension in $Art_{k}$ is a short exact sequence
\begin{equation}
e: 0 \to M \to B \to A \to 0,
\end{equation}
where $B\to A$ is a morphism in $Art_{k}$ and $M$ is an ideal of $B$ annihilated by the maximal ideal $m_{B}$ (so $M$ is a vector space over $k=B/m_{B}$). For simplicity, we also say that $B\to A$ is a small extension.

\end{definition}

If $M=(\eta)$, where $\eta^{2}=0$, then $M$ (as a $k$-vector space) has dimension 1. In this case, the small extension $e$ is called a principal small extension. Since $(\eta) \cong k$, we may also write the principal small extension $e$ as
\begin{equation*}
e: 0 \to k \to B \to A \to 0.
\end{equation*}

A principal small extension is called curvilinear, if it has the form
\begin{equation*}
0 \to (t^{n}) \rightarrow k[t]/(t^{n+1}) \to k[t]/(t^{n}) \to 0,
\end{equation*}
where $n$ is some positive integer.

\begin{definition} [Definition 3.1 of \cite{FM-Obs}] \label{Definition: obs theory}
Let $F$ be a functor of Artin rings. An obstruction theory $(V, v_{e})$ for $F$ is the data of an obstruction space $V \in Set_{*}$ and, for every small extension $e$ (2.1), of an obstruction map $v_{e}: F(A) \times \check{M} \to V$, where $\check{M}$ is the dual space of the $k$-vector space $M$. The obstruction map $v_{e}$ must satisfy the following two conditions:
\begin{itemize}
\item [$\mathrm{(1)}$] \ $v_{e}(*,1)=*$, where $1 \in \check{M}$ is the identity, 
\item [$\mathrm{(2)}$] \ $\mathrm{(base \ change)}$ for every morphism of small extensions $e \to e'$, which is a commutative diagram
\begin{equation}
  \begin{CD}
    e: \  0  @>>> M @>>> B @>>> A @>>> 0  \\
     @. @V \gamma_{1}VV  @V \gamma_{2}VV @V \gamma_{3}VV @. \\
    e': \ 0  @>>> M' @>>> B' @>>> A' @>>> 0,
  \end{CD}
\end{equation}
the following diagram commutes
\[
  \begin{CD}
     F(A)\times \check{M'} @>>> F(A')\times \check{M'} \\
     @VVV  @VV v_{e'}V \\
     F(A)\times  \check{M}@> v_{e}>> V.
  \end{CD}
\]
\end{itemize}

\end{definition}

An obstruction theory $(V, v_{e})$ of $F$ is determined by $v_{e}(a,f)$, where $e$ is any small extension (2.1), $a \in F(A)$ and $f \in \check{M}=Hom_{k}(M,k)$. The morphism $f$ induces a morphism of small extensions $f_{*}: e \to e'$
\[
  \begin{CD}
    e: \  0  @>>> M @>>> B @>>> A @>>> 0  \\
     @. @VfVV  @VVV @VV=V @. \\
   e' : \ 0  @>>> k @>>> B'  @>>> A @>>> 0,
  \end{CD}
\]
where $B'$ is the quotient of $B\oplus k$ by the ideal $\{(m, f(m))| m \in M \}$, see the proof of Lemma 1.1 of \cite{FM-Obs} for details of $B'$.

By base change, there exists a commutative diagram
\[
  \begin{CD}
     F(A)\times \check{k} @>=>> F(A)\times \check{k} \\
     @VVV  @VV v_{e'}V \\
     F(A)\times  \check{M}@> v_{e}>> V.
  \end{CD}
\]
This shows that $v_{e'}(a,g)=v_{e}(a, g \circ f)$, where $a \in F(A)$ and $g \in \check{k}=Hom_{k}(k,k)$. In particular, let $g=id$ be the identity map, one sees that $v_{e'}(a, id)=v_{e}(a, f)$. Noting that $e'$ is a principal small extension, one obtains the following:
\begin{remark} [see Remark 4.3 of \cite{Manetti2005}] \label{Remark: obs by principal}
Every obstruction theory of a functor of Artin rings $F$ is determined by its behavior on principal small extensions. 
\end{remark}

The following proposition explains the name ``obstruction theory".
\begin{proposition} [Proposition 3.3 of \cite{FM-Obs}] \label{Prop: why called obs}
Let $F$ be a functor of Artin rings and $(V, v_{e})$ be an obstruction theory  for $F$. For every small extension $e$ (2.1), let $a \in F(A)$, if $a$ is contained in the image of $F(B)$, then $v_{e}(a, f)=*$ for every $f \in \check{M}$.

\end{proposition}

Two important properties of an obstruction theory are completeness and linearity. 
\begin{definition}  [Definition 4.1 of \cite{FM-Obs}] \label{Definition: complete}
An obstruction theory $(V, v_{e})$ for a functor of Artin rings  $F$ is called complete if for any small extension $e$ (2.1), an element $a \in F(A)$ lifts to $F(B)$ if and only if for every $f \in \check{M}$, $v_{e}(a, f)=*$.
\end{definition}

\begin{definition} [Definition 4.7 of \cite{FM-Obs}] \label{Definition: linear}
An obstruction theory $(V, v_{e})$ for a functor of Artin rings $F$ is called linear if $V$ is a $k$-vector space and for any small extension $e$ (2.1), there is a map $v: F(A) \to V \otimes_{k}M$ such that the induced map 
\[
F(A) \times \check{M} \xrightarrow{(v,1)} (V \otimes_{k}M)\times \check{M} \xrightarrow{\tilde{v}} V,
\]
where $\tilde{v}(x\otimes m,f)=f(m)x$, is the obstruction map $v_{e}$.

\end{definition}

In Lemma \ref{Lemma: LocalHIlbert-CompleteLinear} below,  we explain that the local Hilbert functor has a complete and linear obstruction theory.
\begin{remark} \label{Remark: called Obs}
When $(V, v_{e})$ is a complete and linear obstruction theory, the vector space $V$ is the obstruction space in the usual sense. In this case, for $a \in F(A)$, the image $v(a) \in V \otimes_{k}M$ is often called the obstruction to lifting $a$. See Definition 4.7 of \cite{FM-Obs}.
\end{remark}

A morphism $(V, v_{e}) \to (W, w_{e})$ of two obstruction theories of a functor of Artin rings $F$ is a morphism $f: V \to W$ such that $w_{e} =f \circ v_{e}$ for every small extension $e$.

\begin{definition} [see page 549 of \cite{FM-Obs}]  \label{Definition: Universal obs}
Let $F$ be a functor of Artin rings. An obstruction theory $(O_{F}, ob_{e})$ for $F$ is called universal if for every obstruction theory $(V, v_{e})$ for $F$, there exists a unique morphism $(O_{F}, ob_{e}) \to (V, v_{e})$. In this case, $O_{F}$ is called the universal obstruction space.
\end{definition}

\begin{theorem} [Theorem 3.2 of  \cite{FM-Obs}]
Let $F$ be a functor of Artin rings, there exists a unique universal obstruction theory $(O_{F}, ob_{e})$ for $F$.
\end{theorem}

Fantechi-Manetti constructed the universal obstruction space $O_{F}$ as a quotient of a set modulo the equivalence relation from base change (condition (2) in Definition \ref{Definition: obs theory}). To be precise, we denote by $\hat{O}$ the set
\[
\hat{O}=\bigcup_{e} F(A) \times \check{M},
\]
where $e$ is any small extension $0 \to M \to B \to A \to 0$. For a morphism of small extensions $e \to e'$ (see (2.2) on page 5), there is a set-theoretic map 
\[
F(A) \times \check{M'} \to (F(A) \times \check{M}) \times (F(A') \times \check{M'})
\]
defined as
\[
(a,f) \to (a, f \circ \gamma_{1}) \times (\tilde{\gamma_{3}}(a), f),
\]
where $\tilde{\gamma_{3}}: F(A) \to F(A')$ is the map induced by $\gamma_{3}$. We define an equivalence relation $\backsim$ on $\hat{O}$ to be $(a, f \circ \gamma_{1}) \backsim (\tilde{\gamma_{3}}(a), f)$.

The universal obstruction space $O_{F}$ is defined to be the quotient of the set $\hat{O}$ by the equivalence relation $\backsim$.

It follows from Remark \ref{Remark: obs by principal} that the universal obstruction space $O_{F}$ can be constructed from principal small extensions. Concretely, we consider the set $\tilde{O}=\{ (e, a)| e: 0 \to (\eta) \to B \to A \to 0, a \in F(A) \}$, where $e$ is any principal small extension.
\begin{remark} \label{Remark: Universal Obs-Space}
Let $1_{\eta}: (\eta) \to k$ be the $k$-linear morphism mapping $\eta$ to $1$, the set $\tilde{O}$ can be identified with a subset of $\hat{O}$ by the map $(e,a) \to (e,a, 1_{\eta})$. 

Let $\thickapprox$ be the restriction of $\thicksim$ to $\tilde{O}$, the quotient of $\tilde{O}$ by the equivalence relation $\thickapprox$ is isomorphic to the universal obstruction space $O_{F}$ (see Remark on page 550 of \cite{FM-Obs} for details). This gives another description of universal obstruction space $O_{F}$, which is used below.
\end{remark}

 \begin{definition} \label{Definition: trivial obs}
An obstruction theory $(V, v_{e})$ for a functor of Artin rings $F$ is called trivial, if $V=\{* \}$.
\end{definition}
It is clear that a trivial obstruction theory is linear in the sense of Definition \ref{Definition: linear}.

\begin{definition} \label{Definition: smooth functor}
A functor of Artin rings $F$ is smooth, if $F(B) \to F(A)$ is surjective for every small extension 
\[
e: 0 \to M \to B \to A \to 0.
\]
\end{definition}

\begin{lemma} [Remark on page 552 of \cite{FM-Obs}] \label{Lemma: smooth functor&trivial obs}
If a functor of Artin rings $F$ is smooth, then the universal obstruction theory for $F$ is trivial, which is complete and linear. 
 \end{lemma}

To give an example of obstruction theory, we recall the local Hilbert functor. Let $X$ be a smooth projective variety over a field $k$ of characteristic zero and let $Y \subset X$ be a closed subvariety of codimension $p$, which is a locally complete intersection. For any $A \in Art_{k}$, we denote by $X \times A$ the fibre product $X \otimes_{\mathrm{Spec}(k)} \mathrm{Spec}(A)$. An infinitesimal embedded deformation of $Y$ in $X \times A$ is a closed subscheme $Y^{'} \subset X \times A$ such that $Y^{'}$ is flat over $\mathrm{Spec}(A)$ and $Y^{'} \otimes_{\mathrm{Spec}(A)} \mathrm{Spec}(k)=Y$. 
\begin{definition} \label{Definition: local Hilbert}
The local Hilbert functor $\mathrm{\mathbb{H}ilb}$ is defined to be\footnote{The functor $\mathrm{\mathbb{H}ilb}$ can be defined without assuming that $X$ is smooth projective or $Y$ is a locally complete intersection, for example, see section 3.2 of \cite{Sernesi}.}
\begin{align*}
\mathrm{\mathbb{H}ilb}: \  Art_{k} & \longrightarrow  Set_{*} \\
  \ A & \longrightarrow \mathrm{\mathbb{H}ilb}(A),
\end{align*}
where $\mathrm{\mathbb{H}ilb}(A)$ denotes the set of infinitesimal embedded deformations of $Y$ in $X \times A$.
\end{definition}
The functor $\mathrm{\mathbb{H}ilb}$ is a well-known example of functors of Artin rings, which was intensively studied in literatures. We briefly recall obstruction theory of the functor $\mathrm{\mathbb{H}ilb}$ and refer to \cite{ Hartshorne2, Sernesi, Vistoli} for details.

Since $Y \subset X$ is a locally complete intersection of codimension $p$, there exists a finite open affine covering $\{U_{i} \}_{i \in I}$ of $X$ such that $Y \cap U_{i}$ is given by a regular sequence $f_{i1}, \cdots, f_{ip}$ of $O_{X}(U_{i})$. By Remark \ref{Remark: obs by principal}, it suffices to consider the obstructions arising from principal small extensions. For any principal small extension
\[
e: 0 \to (\eta) \to B \to A \to 0,
\]
and for $Y' \in \mathrm{\mathbb{H}ilb}(A)$, $Y'$ is still a locally complete intersection.  In fact, $Y' \cap U_{i}$ is given by a regular sequence $f^{A}_{i1}, \cdots, f^{A}_{ip}$ of $O_{X}(U_{i})\otimes_{k} A$ (=$O_{X \times A}(U_{i})$).
 
The regular sequence $f^{A}_{i1}, \cdots, f^{A}_{ip}$ can be lifted to a regular sequence $f^{B}_{i1}, \cdots, f^{B}_{ip}$ of $O_{X}(U_{i}) \otimes_{k} B$. On the intersection $U_{ij}=U_{i} \cap U_{j}$, two liftings $f^{B}_{i1}, \cdots, f^{B}_{ip}$ and $f^{B}_{j1}, \cdots, f^{B}_{jp}$ satisfy that $f^{B}_{il}-f^{B}_{jl}=\eta h_{ijl}$, where $h_{ijl} \in O_{X}(U_{ij})$ and $l=1,\cdots,p$. 

Let $\mathcal{I}$ be the ideal sheaf of $Y$ in $X$, the normal bundle $N_{Y/X}$ is defined to be $\mathcal{H}om_{O_{Y}}(\mathcal{I}/\mathcal{I}^{2}, O_{Y})$. The image of $h_{ijl}$ in $O_{Y}(U_{ij})$, denoted $h'_{ijl}$, defines a morphism
\begin{align*}
\mu_{ij}: \  \mathcal{I}(U_{ij})& \to O_{Y}(U_{ij}) \\
f_{il} & \to h'_{ijl}, \mathrm{where} \ l=1,\cdots,p.
\end{align*}
Using $\mathcal{H}om_{O_{X}}(\mathcal{I}, O_{Y})=\mathcal{H}om_{O_{Y}}(\mathcal{I}/\mathcal{I}^{2}, O_{Y})$, we consider $\mu_{ij}$ as an element of $\Gamma(U_{ij}, N_{Y/X})$, which forms a \v{C}ech 1-cocycle $(\mu_{ij})_{i,j}$.
It is known that the 1-cocycle $(\mu_{ij})_{i,j}$ gives a well-defined class\footnote{It does not depends on the choice of lifting of $f^{A}_{i1}, \cdots, f^{A}_{ip}$.} of $H^{1}(X, N_{Y/X})$, denoted $\{\mu_{ij} \}_{i,j}$, which is the obstruction to lifting $Y'$.

We use the class $\{\mu_{ij} \}_{i,j}$ to define a set-theoretic map
\begin{align*}
v: \mathrm{\mathbb{H}ilb}(A) & \to H^{1}(X, N_{X/Y})\otimes_{k}(\eta) \\
Y' & \longrightarrow \{\mu_{ij} \}_{i,j} \otimes \eta.
\end{align*}
The map $v$ induces a set-theoretic map
\begin{equation}
v_{e}: \mathrm{\mathbb{H}ilb}(A) \times \check{(\eta)}  \to H^{1}(X, N_{X/Y}),
\end{equation}
which is defined to be the composition
\[
\mathrm{\mathbb{H}ilb}(A) \times \check{(\eta)}  \xrightarrow{(v,1)} (H^{1}(X, N_{X/Y})\otimes_{k}(\eta)) \times \check{(\eta)} \xrightarrow{\tilde{v}} H^{1}(X, N_{X/Y}),
\]
where $\check{(\eta)}$ is the dual space of the $k$-vector space $(\eta)$ and for $g \in \check{(\eta)}$ and $x \in H^{1}(X, N_{X/Y})$, $\tilde{v}$ is defined by $\tilde{v}(x \otimes \eta, g)=g(\eta)x$. Concretely,
\begin{equation}
v_{e}(Y',g)=g(\eta)\{\mu_{ij} \}_{i,j}.
\end{equation}

Now, we show that $v_{e}$ is an obstruction map in the sense of Definition \ref{Definition: obs theory}. It suffices to check that $v_{e}$ satisfies base change. For every morphism of principal small extensions $e \to e'$ 
\[
  \begin{CD}
    e: \  0  @>>> (\eta) @>>> B @>>> A @>>> 0  \\
     @. @V\gamma_{1}VV  @V\gamma_{2}VV @V\gamma_{3}VV @. \\
    e': \ 0  @>>> (\eta') @>>> B' @>>> A' @>>> 0,
  \end{CD}
\]
the morphisms $\gamma_{3}$ and $\gamma_{2}$ induce maps $\tilde{\gamma_{3}}: \mathrm{\mathbb{H}ilb}(A) \to \mathrm{\mathbb{H}ilb}(A')$ and $\tilde{\gamma_{2}}: \mathrm{\mathbb{H}ilb}(B) \to\mathrm{\mathbb{H}ilb}(B')$ respectively.

For $Y' \in \mathrm{\mathbb{H}ilb}(A)$ which is locally on $U_{i}$ given by $f^{A}_{i1}, \cdots, f^{A}_{ip}$, $\tilde{\gamma_{3}}(Y') \in \mathrm{\mathbb{H}ilb}(A')$ is locally on $U_{i}$ given by the sequence $1\otimes \gamma_{3}(f^{A}_{i1}), \cdots, 1\otimes \gamma_{3}(f^{A}_{ip})$, where $1\otimes \gamma_{3}: O_{X}(U_{i})\otimes_{k} A \to O_{X}(U_{i})\otimes_{k} A'$ is the morphism induced by $\gamma_{3}$. Since the sequence $f^{B}_{i1}, \cdots, f^{B}_{ip}$ is a lifting of the sequence $f^{A}_{i1}, \cdots, f^{A}_{ip}$, the sequence $1\otimes \gamma_{2}(f^{B}_{i1}), \cdots, 1\otimes\gamma_{2}(f^{B}_{ip})$ is a lifting of the sequence $1\otimes \gamma_{3}(f^{A}_{i1}), \cdots, 1\otimes \gamma_{3}(f^{A}_{ip})$, where $1\otimes \gamma_{2}: O_{X}(U_{i})\otimes_{k} B \to O_{X}(U_{i})\otimes_{k} B'$ is the morpism induced by $\gamma_{2}$. 

On the intersection $U_{ij}$, the difference of two liftings $1\otimes \gamma_{2}(f^{B}_{i1}), \cdots, 1\otimes \gamma_{2}(f^{B}_{ip})$ and $1\otimes \gamma_{2}(f^{B}_{j1}), \cdots, 1\otimes \gamma_{2}(f^{B}_{jp})$ satisfies that, for $l=1,\cdots,p$,
\begin{align*}
1\otimes \gamma_{2}(f^{B}_{il})-1\otimes\gamma_{2}(f^{B}_{jl})& =1\otimes \gamma_{2}(f^{B}_{il}- f^{B}_{jl})=1 \otimes \gamma_{2}(h_{ijl} \eta) \\
& =h_{ijl}\gamma_{2}(\eta)=h_{ijl}\gamma_{1}(\eta).
\end{align*}

For $g \in \check{(\eta')}$, $g \circ \gamma_{1} \in \check{(\eta)}$, it follows from the definition of $v_{e}$ (see (2.3) and (2.4)) that $v_{e'}(\tilde{\gamma_{3}}(Y'), g)=g(\gamma_{1}(\eta))\{\mu_{ij}\}_{i,j}$ and $v_{e}(Y', g \circ \gamma_{1})=g \circ \gamma_{1}(\eta)\{\mu_{ij}\}_{i,j}=g(\gamma_{1}(\eta))\{\mu_{ij}\}_{i,j}$. This implies a commutative diagram
\[
  \begin{CD}
     (Y', g) @>(\tilde{\gamma_{3}}, 1)>>  (\tilde{\gamma_{3}}(Y'), g)\\
     @V(1, \tilde{\gamma_{1}})VV  @V v_{e'}VV \\
      (Y', g \circ \gamma_{1})@>v_{e}>> v_{e'}(\tilde{\gamma_{3}}(Y'), g),
  \end{CD}
\]
where $\tilde{\gamma_{1}}: \check{(\eta')} \to \check{(\eta)}$ is the morphism induced by $\gamma_{1}$. It follows that the map $v_{e}$ (2.3) satisfies base change and it is an obstruction map. Hence,  $(H^{1}(X, N_{X/Y}), v_{e})$ is an obstruction theory of the local Hilbert functor $\mathrm{\mathbb{H}ilb}$.

It can be easily checked that the obstruction theory $(H^{1}(X, N_{X/Y}), v_{e})$ is complete and linear (in the sense of Definition \ref{Definition: complete} and Definition \ref{Definition: linear} respectively). In summary,
\begin{lemma}\label{Lemma: LocalHIlbert-CompleteLinear}
With notation as above, the local Hilbert functor $\mathrm{\mathbb{H}ilb}$ has a complete and linear obstruction theory $(H^{1}(X, N_{X/Y}), v_{e})$.
\end{lemma}

By Remark \ref{Remark: Universal Obs-Space}, an element of the universal obstruction space $O_{\mathrm{\mathbb{H}ilb}}$ of the local Hilbert functor $\mathrm{\mathbb{H}ilb}$ is an equivalent class $[(e, Y')]$, where $e$ is a principal small extension $0 \to (\eta) \to B \to A \to 0$ and $Y' \in \mathrm{\mathbb{H}ilb}(A)$. The obstruction map $v_{e}$ (2.3) induces a set-theoretic map
\begin{align}
[v_{e}]: O_{\mathrm{\mathbb{H}ilb}} & \to H^{1}(Y, N_{Y/X}) \\
[(e, Y')] & \to v_{e}(Y', 1_{\eta}),   \notag
\end{align}
where $1_{\eta}: (\eta) \to k$ is the morphism mapping $\eta$ to 1 (see Remark \ref{Remark: Universal Obs-Space}). It follows from (2.4) that $v_{e}(Y', 1_{\eta})=1_{\eta}(\eta)\{\mu_{ij}\}_{i,j}=1\{\mu_{ij}\}_{i,j}=\{\mu_{ij}\}_{i,j}$.

Fibred products exist on the category $Art_{k}$. Given morphisms $C \xrightarrow{f} A$ and $B \xrightarrow{g} A$ in $Art_{k}$, there exists a commutative diagram
\[
  \begin{CD}
     B \times_{A} C @>>> C \\
     @VVV  @VVfV \\
     B @>g>> A,
  \end{CD}
\]
where $ B \times_{A} C= \{(b,c)| g(b)=f(c) \} $ is also in the category $Art_{k}$.

Let $F$ be a functor of Artin rings, there is a natural map
\begin{equation}
S: F(B \times_{A} C) \to F(B) \times_{F(A)} F(C).
\end{equation}

\begin{definition} [see Definition 2.7 of \cite{FM-Obs}]\label{Definition: Schelessinger cond}
Let $F$ be a functor of Artin rings, the followings are called Schlessinger conditions:
\begin{itemize}
\item [$\mathrm{(H1)}$] \  Map $S$ is surjective if $C \to A$ is a small extension. 

\item [$\mathrm{(H2)}$] \  Map $S$ is bijective if $A=k$, and $C=k[\varepsilon]$ is the ring of dual numbers. 

\item [$\mathrm{(H3)}$] \ Conditions $\mathrm{(H1)}$ and $\mathrm{(H2)}$ hold and the dimension of $F(k[\varepsilon])$ is finite\footnote{The condition (H2) guarantees that the set $F(k[\varepsilon])$ carries a structure of $k$-vector space, so it makes sense to speak about the dimension of $F(k[\varepsilon])$ in (H3).}.

\item [$\mathrm{(H4)}$] \  Map $S$ is bijective if $C \to A$ is a small extension.
\end{itemize}
\end{definition}

The following theorem shows that, for certain functors of Artin rings, the universal obstruction spaces carry structures of vector spaces. 

\begin{theorem} [Theorem 6.6 of \cite{FM-Obs}] \label{Theorem: obs space is VS}
Assume that a functor of Artin rings $F$, which satisfies the Schlessinger conditions $\mathrm{(H1)}$ and $\mathrm{(H2)}$\footnote{Such a functor is called a functor with good deformation theory, see Definition 2.8 in \cite{FM-Obs} and references therein.}, has a complete and linear obstruction theory $(V, v_{e})$. Then there exists a unique structure of vector space on $O_{F}$ such that the universal obstruction theory $(O_{F}, ob_{e})$ is linear in the sense of Definition \ref{Definition: linear}. Moreover, the induced map $O_{F} \to V$ is a linear monomorphism of vector spaces.
\end{theorem}

It is well known that the local Hilbert functor $\mathrm{\mathbb{H}ilb}$ in Definition \ref{Definition: local Hilbert} satisfies the Schlessinger conditions $\mathrm{(H1)}$, $\mathrm{(H2)}$, $\mathrm{(H3)}$ and $\mathrm{(H4)}$.
Combing Lemma \ref{Lemma: LocalHIlbert-CompleteLinear} with Theorem \ref{Theorem: obs space is VS}, one deduces that
\begin{corollary} \label{Corollary: Obs HIlb is VS}
The universal obstruction space $O_{\mathrm{\mathbb{H}ilb}}$ of the local HIlbert functor $\mathrm{\mathbb{H}ilb}$ is a $k$-vector space and the set-theoretic map (2.5)
\[
[v_{e}]: O_{\mathrm{\mathbb{H}ilb}} \rightarrow H^{1}(Y,N_{Y/X}),
\]
is a linear monomorphism of $k$-vector spaces.
\end{corollary}
 
From the definition of the map $[v_{e}]$ (2.5), one sees that the image of $[v_{e}]$ is obstructions to embedded deformations of $Y$ in $X$. Hence, an alternative way to state Conjecture \ref{Conjecture:main} is as following.
\begin{conjecture} \label{Conjecture:main-2}
The composition of morphisms of $k$-vector spaces
 \[
O_{\mathrm{\mathbb{H}ilb}} \xrightarrow{[v_{e}]} H^{1}(N_{Y/X}) \to H^{p+1}(\Omega_{X/k}^{p-1})
 \]
 is trivial, where $H^{1}(N_{Y/X}) \to H^{p+1}(\Omega_{X/k}^{p-1})$ is the semi-regularity map.
\end{conjecture}

\section{Chern character and Hochschild homology}
In this section, after recalling background on Hochschild homology in section 3.1, we compute Hochschild homology in section 3.2. Then we introduce a functor $\mathbb{HH}^{(p)}$ in Definition \ref{Definition: Functor-HH-final} and construct a natural transformation from the local Hilbert functor to the functor $\mathbb{HH}^{(p)}$ in Lemma \ref{Lemma: transf-Hilb-HH-final}.

\subsection{K-theory and Hochschild homology}

Let $X$ be a noetherian scheme of finite type over a field $k$ of characteristic zero, the Hochschild homology complexes  $HH(X)$ and $HH^{\mathbb{Q}}(X)$ (over $\mathbb{Q}$) and the negative cyclic homology complexes $HN^{\mathbb{Q}}(X)$ (over $\mathbb{Q}$) are defined from a localization pair. Let $Y \subset X$ be closed, the Hochschild homology complex $HH(X \ \mathrm{on} \  Y )$ and the negative cyclic homology complex $HN(X \ \mathrm{on} \  Y )$ can be defined in a similar way. We refer to Example 2.7 and 2.8 of \cite{CHSW} for details.

Following the convention in section 2 of \cite{CHSW}, we use cohomological notation to define the Hochschild homology with support $HH_{i}(X \ \mathrm{on} \  Y )$ to be
\begin{equation}
HH_{i}(X \ \mathrm{on} \  Y ):=H^{-i}(HH(X \ \mathrm{on} \  Y )).
\end{equation}

\begin{definition} \label{Definition: Hochschild-Sheaf}
With notation as above, one defines $\mathcal{HH}_{0}(O_{X} \ \mathrm{on} \ Y)$ to be the sheaf associated to the presheaf
\[
U \to HH_{0}(U \ \mathrm{on} \  Y \cap U),
\]
where $U \subseteq X$ is open and $HH_{0}(U \ \mathrm{on} \  Y \cap U)$ is defined in (3.1).
\end{definition}

\begin{definition} \label{Definition:K-Sheaf}
 With notation as above, one defines $\mathcal{K}_{0}(O_{X} \ \mathrm{on} \ Y)$ to be the sheaf associated to the presheaf
\[
U \to K_{0}(U \ \mathrm{on} \  Y \cap U),
\]
where $U \subseteq X$ is open and $K_{0}(U \ \mathrm{on} \  Y \cap U)$ is Grothendieck group of derived category obtained from the exact category of perfect complexes of $O_{U}$-modules supported on $Y \cap U$.
\end{definition}

Let $\mathcal{HH}(X)$, $\mathcal{HH}^{\mathbb{Q}}(X)$ and $\mathcal{HN}^{\mathbb{Q}}(X)$ be the Eilenberg-Mac Lane spectra associated to the complexes $HH(X)$, $HH^{\mathbb{Q}}(X)$ and $HN^{\mathbb{Q}}(X)$ respectively. According to Corti\~nas, Haesemeyer, Schlichting and Weibel \cite{CHSW} (page 565), there exists a Chern character $\mathcal{K}(X) \to \mathcal{HN}^{\mathbb{Q}}(X)$, where $\mathcal{K}(X)$ is the non-connective K-theory spectrum of perfect complexes on $X$. Composing $\mathcal{K}(X) \to \mathcal{HN}^{\mathbb{Q}}(X)$ with the maps $\mathcal{HN}^{\mathbb{Q}}(X) \to \mathcal{HH}^{\mathbb{Q}}(X) \to \mathcal{HH}(X)$, one obtains a map of spectra $\mathcal{K}(X) \to \mathcal{HH}(X)$, which induces a map
\begin{equation}
\mathcal{K}(X \ \mathrm{on} \ Y) \to \mathcal{HH}(X \ \mathrm{on} \ Y),
\end{equation}
where $\mathcal{K}(X \ \mathrm{on} \ Y)$ is the non-connective K-theory spectrum of perfect complexes of $O_{X}$-modules supported on $Y$ and $\mathcal{HH}(X \ \mathrm{on} \ Y)$ is the Eilenberg-Mac Lane spectrum associated to the Hochschild homology complex $HH(X \ \mathrm{on} \ Y)$.

From now on, $X$ is a smooth projective variety over a field $k$ of characteristic zero and  $Y \subset X$ is a closed  subvariety of codimension $p$, which is a locally complete intersection. The following setting is frequently used below.
\begin{setting} \label{Setting:S}
Let $\{U_{i} \}_{i \in I}$ be an open affine covering of $X$ such that $Y \cap U_{i}$ is defined by a regular sequence $f_{i1}, \cdots, f_{ip}$. For any $A \in Art_{k}$ and for $Y' \in \mathrm{\mathbb{H}ilb}(A)$,  $Y' \cap U_{i}$ is also given by a a regular sequence, denoted $f^{A}_{i1}, \cdots, f^{A}_{ip}$.

For any $A \in Art_{k}$, we abbreviate $X \otimes_{\mathrm{Spec}(k)} \mathrm{Spec}(A)$ to $X \times A$.

\end{setting}

The element $Y' \in \mathrm{\mathbb{H}ilb}(A)$ gives an element of $H^{0}(X, \mathcal{K}_{0}(O_{X \times A} \ \mathrm{on} \ Y))$ as follows. Let $L^{i,A}_{\bullet}$ be the Koszul complex of the regular sequence $f^{A}_{i1}, \cdots, f^{A}_{ip}$, we consider $L^{i,A}_{\bullet}$ as an element of $\Gamma(U_{i}, \mathcal{K}_{0}(O_{X \times A} \  \mathrm{on} \ Y))$. The restriction of $L^{i,A}_{\bullet}$ on the intersection $U_{ij}=U_{i} \cap U_{j}$, denoted $L^{i,A}_{\bullet} |_{U_{ij}}$, satisfies that
\[
L^{i,A}_{\bullet}|_{U_{ij}}=L^{j,A}_{\bullet}|_{U_{ij}} \in \Gamma(U_{i} \cap U_{j}, \mathcal{K}_{0}(O_{X \times A} \ \mathrm{on} \ Y)),
\]
since both of them are the resolution of $Y' \cap U_{ij}$. Hence, $(L^{i,A}_{\bullet})_{i}$ is a cocycle and then gives an element of $H^{0}(X, \mathcal{K}_{0}(O_{X \times A} \ \mathrm{on} \ Y))$.

\begin{definition} \label{Definition:Hilb-Map-K}
With notation as above, one defines a set-theoretic map 
\begin{align}
\alpha_{A}: \  \mathrm{\mathbb{H}ilb}(A) &  \longrightarrow H^{0}(X, \mathcal{K}_{0}(O_{X \times A} \ \mathrm{on} \ Y)) \\
Y^{'} \ &  \longrightarrow   (L^{i,A}_{\bullet})_{i}.    \notag
\end{align}

\end{definition}

The map (3.2)
\[
\mathcal{K}(X \times A \ \mathrm{on} \ Y) \to \mathcal{HH}(X \times A \ \mathrm{on} \ Y)
\]
induces a map of groups
\begin{equation}
ch_{A}: H^{0}(X, \mathcal{K}_{0}(O_{X \times A} \ \mathrm{on} \ Y)) \to H^{0}(X, \mathcal{HH}_{0}(O_{X \times A} \ \mathrm{on} \ Y)).
\end{equation}

\begin{definition} \label{Definition: Relative Hochschild}
One defines $\mathcal{\overline{HH}}_{0}(O_{X \times A} \ \mathrm{on} \  Y)$ to be the kernel of the morphism (induced by the augmentation map $A \to k$)
\[
 \mathcal{HH}_{0}(O_{X \times A} \ \mathrm{on} \ Y) \to \mathcal{HH}_{0}(O_{X} \ \mathrm{on} \ Y).
\]
\end{definition}

Since the morphism $\mathcal{HH}_{0}(O_{X \times A} \ \mathrm{on} \ Y) \to \mathcal{HH}_{0}(O_{X} \ \mathrm{on} \ Y)$ is split by $\mathcal{HH}_{0}(O_{X} \ \mathrm{on} \ Y) \to \mathcal{HH}_{0}(O_{X \times A} \ \mathrm{on} \ Y)$, it follows that
\[
 \mathcal{HH}_{0}(O_{X \times A} \ \mathrm{on} \ Y) =  \mathcal{HH}_{0}(O_{X} \ \mathrm{on} \ Y) \oplus \mathcal{\overline{HH}}_{0}(O_{X \times A} \ \mathrm{on} \  Y).
\]
Hence, there is a natural projection map
\begin{equation}
r_{A}: H^{0}(X, \mathcal{HH}_{0}(O_{X \times A} \ \mathrm{on} \ Y)) \to H^{0}(X, \mathcal{\overline{HH}}_{0}(O_{X \times A} \ \mathrm{on} \  Y)).
\end{equation}

\begin{definition} \label{Definition: Hilb-Map-Relative HH}
With notation as above, one defines a set-theoretic map 
\[
\mathrm{\mathbb{H}ilb}(A) \to H^{0}(X, \mathcal{\overline{HH}}_{0}(O_{X \times A} \ \mathrm{on} \  Y))
\]
to be the composition
{\footnotesize
\begin{align*}
&\mathrm{\mathbb{H}ilb}(A) \xrightarrow{\alpha_{A}(3.3)} H^{0}(X, \mathcal{K}_{0}(O_{X \times A} \ \mathrm{on} \ Y)) \xrightarrow{ch_{A}(3.4)} H^{0}(X, \mathcal{HH}_{0}(O_{X \times A} \ \mathrm{on} \ Y))  \\
& \xrightarrow{r_{A}(3.5)} H^{0}(X, \mathcal{\overline{HH}}_{0}(O_{X \times A} \ \mathrm{on} \  Y)). 
\end{align*}
}
\end{definition}

Let $f: B \to A$ be a morphism in the category $Art_{k}$, there exists a commutative diagram of sets
\begin{equation}
\begin{CD}
\mathrm{\mathbb{H}ilb}(B) @>\alpha_{B}>>  H^{0}(X,  \mathcal{K}_{0}(O_{X \times B} \ \mathrm{on} \  Y))  \\
@Vf_{Hilb}VV @Vf_{K_{0}}VV \\
\mathrm{\mathbb{H}ilb}(A) @>\alpha_{A}>>  H^{0}(X,  \mathcal{K}_{0}(O_{X \times A} \ \mathrm{on} \  Y)), \\
\end{CD}
\end{equation}
where $f_{Hilb}$ and $f_{K_{0}}$ are induced by $f$ respectively. To check the commutativity, for $Y' \in \mathrm{\mathbb{H}ilb}(B)$, in notation of Definition \ref{Definition:Hilb-Map-K}, one sees that $\alpha_{B}(Y')=(L^{i,B}_{\bullet})_{i}$ and 
\[
f_{K_{0}} \circ \alpha_{B}(Y')=(f^{*}L^{i,B}_{\bullet})_{i}=(L^{i,B}_{\bullet} \otimes_{B}A)_{i},
\]
where $f^{*}$ is the pull-back: {\small $K_{0}(O_{X \times B}(U_{i}) \ \mathrm{on} \  Y \cap U_{i}) \to K_{0}(O_{X \times A}(U_{i}) \ \mathrm{on} \  Y \cap U_{i})$}, see (3.14) on page 317 of \cite{TT}.

Since $f_{Hilb}(Y')=Y' \times_{\mathrm{Spec}(B)}\mathrm{Spec}(A) \in \mathrm{\mathbb{H}ilb}(A)$, the Koszul resolution of $f_{Hilb}(Y') \cap U_{i}$ is the Koszul complex $L^{i,B}_{\bullet} \otimes_{B}A$. Hence, there is a commutative diagram
\[
\begin{CD}
Y' @>\alpha_{B}>> (L^{i,B}_{\bullet})_{i} \\
@Vf_{Hilb}VV @Vf_{K_{0}}VV   \\
Y' \times_{\mathrm{Spec}(B)}\mathrm{Spec}(A) @>\alpha_{A}>>  (L^{i,B}_{\bullet} \otimes_{B}A)_{i},
\end{CD}
\]
which shows the commutativity of (3.6).

There also exists a commutative diagram
{\tiny
\begin{equation}
\begin{CD}
 H^{0}( \mathcal{K}_{0}(O_{X \times B} \ \mathrm{on} \  Y)) @>ch_{B}>>  H^{0}(\mathcal{HH}_{0}(O_{X \times B} \ \mathrm{on} \ Y))  @>r_{B}>> H^{0}(\mathcal{\overline{HH}}_{0}(O_{X \times B} \ \mathrm{on} \  Y)) \\
 @Vf_{K_{0}}VV  @Vf_{HH}VV  @V\overline{f}_{HH}VV \\
 H^{0}(\mathcal{K}_{0}(O_{X \times A} \ \mathrm{on} \  Y)) @>ch_{A}>>  H^{0}(\mathcal{HH}_{0}(O_{X \times A} \ \mathrm{on} \ Y))  @>r_{A}>> H^{0}(\mathcal{\overline{HH}}_{0}(O_{X \times A} \ \mathrm{on} \  Y)), \\
\end{CD}
\end{equation}
}where $f_{HH}$ and $\overline{f}_{HH}$ are the maps induced by $f$ respectively. The commutativity of the left square follows from the naturality of the Chern character and the commutativity of the right square is obvious.

\subsection{Computation of Hochschild homology}

For $S$ a commutative algebra over a field $k$ of characteristic zero, lambda operations $\lambda^{m}$ and Adams operations $\psi^{m}$ are defined for Hochschild homology $HH_{*}(S)$, we refer to section 4.5 of \cite{L-2} and section 9.4.3 of \cite{W1} for details. For $l$ a non-negative integer, let $HH^{(l)}_{*}(S)$ denote the eigenspace of $\psi^{m}=m^{l+1}$.

\begin{lemma} [see Theorem 4.5.12 of \cite{L-2} or Ex 9.4.4 of \cite{W1}]\label{lemma: l-l-omega}
With notation as above, there is an isomorphism $HH^{(l)}_{l}(S) \cong \Omega^{l}_{S/k}$.
\end{lemma}

We keep the notation of Setting \ref{Setting:S}. Let $U_{i}=\mathrm{Spec}(R) \subset X$ be open affine such that $Y \cap U_{i}$ is given by a regular sequence $f_{i1}, \cdots, f_{ip}$, which generates an ideal $J \subset R$.

Let $HH(R \otimes_{k} A)$ be the Hochschild homology complex of $R \otimes_{k} A$. Since Hochschild homology satisfies Zariski descent, the group $HH_{0}(R \otimes_{k} A \ \mathrm{on} \  J)$ defined in (3.1) can be identified with the hypercohomogy of the complex $HH(R \otimes_{k} A)$
\begin{equation*}
HH_{0}(R \otimes_{k} A \ \mathrm{on} \  J)=\mathbb{H}^{0}_{J }(R, HH(R \otimes_{k} A)).
\end{equation*}
This enables us to extend Adams operations $\psi^{m}$ from Hochschild homology of $R \otimes_{k} A$ to $HH_{0}(R \otimes_{k} A \ \mathrm{on} \  J)$. In fact, the complex $HH(R \otimes_{k} A)$ naturally splits into sum of sub-complexes $HH^{(l)}(R \otimes_{k} A)$. The eigenspace (of $HH_{0}(R \otimes_{k} A \ \mathrm{on} \  J)$) of $\psi^{m}=m^{l+1}$, denoted $HH^{(l)}_{0}(R \otimes_{k} A \ \mathrm{on} \  J)$, is the hypercohomogy $\mathbb{H}^{0}_{J }(R, HH^{(l)}(R \otimes_{k} A))$.

It is obvious that Adams operations $\psi^{m}$ can be defined for the sheaf $\mathcal{HH}_{0}(O_{X \times A} \ \mathrm{on} \ Y)$ and we denote by $\mathcal{HH}^{(l)}_{0}(O_{X \times A} \ \mathrm{on} \ Y)$ the eigenspace of $\psi^{m}=m^{l+1}$.

It is noted that $X$ and $X \times A$ have the same underlying space.
\begin{theorem} \label{Theorem: Compute-HH-2}
With notation as above, there is an isomorphism of sheaves
\begin{align*}
\mathcal{HH}^{(l)}_{0}(O_{X \times A} \ \mathrm{on} \  Y ) = \bigoplus_{ \tiny \begin{array}{c}
j_{1}+j_{2}=p \\
j_{1}+l_{2}=l
\end{array}} \mathcal{H}^{p}_{Y}(\Omega^{j_{1}}_{X/k}) \otimes_{k} HH^{(l_{2})}_{j_{2}}(A),
\end{align*}
where $\mathcal{H}^{p}_{Y}( \Omega^{j_{1}}_{X/k})$ is the local cohomology sheaf and $HH^{(l_{2})}_{j_{2}}(A)$ is the eigenspace of $\psi^{m}=m^{l_{2}+1}$ with $\psi^{m}$ Adams operations on Hochschild homology $HH_{j_{2}}(A)$.
\end{theorem}

It suffices to prove the isomorphism locally on $U_{i}=\mathrm{Spec}(R)$. We need to prove that
\begin{align*}
HH^{(l)}_{0}(R \otimes_{k} A \ \mathrm{on} \  J) & = \mathbb{H}^{0}_{J }(R, HH^{(l)}(R \otimes_{k} A)) \\
&= \bigoplus_{ \tiny \begin{array}{c}
j_{1}+j_{2}=p \\
j_{1}+l_{2}=l
\end{array}} H^{p}_{J}(R,  \Omega^{j_{1}}_{R/k}) \otimes_{k} HH^{(l_{2})}_{j_{2}}(A).
\end{align*}

\begin{proof}
There exists a spectral sequence
\[
E_{2}^{i, j} = H^{i}_{J }(R, H^{j}(HH^{(l)}(R \otimes_{k} A))) \Longrightarrow \mathbb{H}^{0}_{J }(R, HH^{(l)}(R \otimes_{k} A)),
\]
where $i+j=0$. Since we use cohomological notation (see (3.1)), $H^{j}(HH^{(l)}(R \otimes_{k} A))=HH^{(l)}_{-j}(R \otimes_{k} A)$, the above spectral sequence can be rewritten as 
\begin{equation}
E_{2}^{i, j} = H^{i}_{J }(R, HH^{(l)}_{j}(R \otimes_{k} A)) \Longrightarrow \mathbb{H}^{0}_{J }(R, HH^{(l)}(R \otimes_{k} A)),
\end{equation}
where the non-negative integers $i$ and $j$ are equal, i.e., $i=j$.

It is known that $R$ can be written as a direct limit of smooth $k$-algebras $\tilde{R}$'s. The Hochschild-Konstant-Rosenberg theorem says that $HH_{j_{1}}(\tilde{R}) \cong \Omega^{j_{1}}_{\tilde{R}/k}$. Both the Hochschild homology and $\Omega^{j_{1}}_{\tilde{R}/k}$ commutes with limit, so it follows that $HH_{j_{1}}(R) \cong \Omega^{j_{1}}_{R/k}$. On the other hand, it follows from Lemma \ref{lemma: l-l-omega} that $HH^{(j_{1})}_{j_{1}}(R) \cong \Omega^{j_{1}}_{R/k}$. This implies that $HH^{(i)}_{j_{1}}(R) =0$, if $ i \neq j_{1}$.

By K$\ddot{\mathrm{u}}$nneth formula for Hochschild homology in \cite{Kassel}, one has
\begin{align*}
HH^{(l)}_{j}(R \otimes_{k} A) & = \bigoplus_{ \tiny \begin{array}{c}
j_{1}+j_{2}=j \\
l_{1}+l_{2}=l
\end{array}}  HH^{(l_{1})}_{j_{1}}(R) \otimes_{k} HH^{(l_{2})}_{j_{2}}(A) \\
& = \bigoplus_{ \tiny \begin{array}{c}
j_{1}+j_{2}=j \\
j_{1}+l_{2}=l
\end{array}}  HH^{(j_{1})}_{j_{1}}(R) \otimes_{k} HH^{(l_{2})}_{j_{2}}(A) \\
& = \bigoplus_{ \tiny \begin{array}{c}
j_{1}+j_{2}=j \\
j_{1}+l_{2}=l
\end{array}} \Omega^{j_{1}}_{R/k} \otimes_{k} HH^{(l_{2})}_{j_{2}}(A).
\end{align*}

Passing to local cohomology,
\begin{align*}
H^{i}_{J}(R, HH^{(l)}_{j}(R \otimes_{k} A)) & =  \bigoplus_{ \tiny \begin{array}{c}
j_{1}+j_{2}=j \\
j_{1}+l_{2}=l
\end{array}} H^{i}_{J}(R, \Omega^{j_{1}}_{R/k} \otimes_{k} HH^{(l_{2})}_{j_{2}}(A)) \\
& = \bigoplus_{ \tiny \begin{array}{c}
j_{1}+j_{2}=j \\
j_{1}+l_{2}=l
\end{array}} H^{i}_{J}(R, \Omega^{j_{1}}_{R/k}) \otimes_{k} HH^{(l_{2})}_{j_{2}}(A).
\end{align*}

Since the $J$-depth of $\Omega^{j_{1}}_{R/k}$ is $p$, $H^{i}_{J }(R,  \Omega^{j_{1}}_{R/k})=0$ unless $i=p$. There is only one non-zero term $H^{p}_{J }(R, HH^{(l)}_{p}(R \otimes_{k} A))$ in the spectral sequence (3.8). This shows that
{\footnotesize
\begin{align}
\mathbb{H}^{0}_{J }(R, HH^{(l)}(R \otimes A)) & =H^{p}_{J }(R, HH^{(l)}_{p}(R \otimes_{k} A)) \\
& = \bigoplus_{ \tiny \begin{array}{c}   \notag
j_{1}+j_{2}=p \\
j_{1}+l_{2}=l
\end{array}} H^{p}_{J }(R,  \Omega^{j_{1}}_{R/k}) \otimes_{k} HH^{(l_{2})}_{j_{2}}(A).
\end{align}
}
\end{proof}

When $l=p$, then $l_{2}= l-j_{1}= p-j_{1}=j_{2}$. It follows from Lemma \ref{lemma: l-l-omega} that $HH^{(l_{2})}_{j_{2}}(A)=HH^{(j_{2})}_{j_{2}}(A)= \Omega^{j_{2}}_{A/k}$. This shows that

\begin{theorem} \label{Theorem: Compute-HH-2-p case}
Let $l=p$ in Theorem \ref{Theorem: Compute-HH-2}, there is an isomorphism of sheaves
\begin{align*}
 \mathcal{HH}^{(p)}_{0}(O_{X \times A} \ \mathrm{on} \  Y )  =  \bigoplus_{j_{1}+j_{2}=p }\mathcal{H}^{p}_{Y}( \Omega^{j_{1}}_{X/k}) \otimes_{k} \Omega^{j_{2}}_{A/k}. \\
\end{align*}
\end{theorem}

\begin{remark}  \label{Remark: HH(A)-is-omegaXA}
There is an isomorphism of sheaves
\begin{equation}
\Omega^{p}_{X \times A/k} = \bigoplus_{j_{1}+j_{2}=p}\Omega^{j_{1}}_{X/ k} \otimes_{k} \Omega^{j_{2}}_{A/ k},
\end{equation}
which can be checked locally. This implies that there is an isomorphism of local cohomology sheaves
\[
\mathcal{H}_{Y}^{p}(\Omega^{p}_{X \times A/ k}) = \bigoplus_{j_{1}+j_{2}=p }\mathcal{H}^{p}_{Y}(\Omega^{j_{1}}_{X/k}) \otimes_{k} \Omega^{j_{2}}_{A/k}.
\]
By Theorem \ref{Theorem: Compute-HH-2-p case}, there is an isomorphism of sheaves\footnote{This isomorphism can be checked alternatively by using Lemma \ref{lemma: l-l-omega}. For $U_{i}=\mathrm{Spec}(R)$, let $l=p$ in (3.9), it can be deduced that $\mathbb{H}^{0}_{J }(R, HH^{(p)}(R \otimes A)) =H^{p}_{J }(R, HH^{(p)}_{p}(R \otimes_{k} A))=H^{p}_{J }(R, \Omega^{p}_{R \otimes A/k})$.}
\[
\mathcal{HH}^{(p)}_{0}(O_{X \times A} \ \mathrm{on} \  Y ) = \mathcal{H}_{Y}^{p}(\Omega^{p}_{X \times A/ k}).
\]
\end{remark}

\subsection{Chern character}

By Definition \ref{Definition:Hilb-Map-K}, $\alpha_{A}(Y')=(L^{i,A}_{\bullet})_{i}$, where $L^{i,A}_{\bullet}$ is Koszul complex of the regular sequence $f^{A}_{i1}, \cdots, f^{A}_{ip}$. Adams operations $\psi^{m}$ for K-theory of perfect complexes defined in \cite{GilletSoule} has the following property.
\begin{lemma} [Prop 4.12  of \cite{GilletSoule}]  \label{Lemma: GilletSoule}
Adams operations $\psi^{m}$ on $L^{i,A}_{\bullet}$ satisfies that
\[
\psi^{m}(L^{i,A}_{\bullet})  = m^{p}L^{i,A}_{\bullet}.
\]
\end{lemma}

Let $\mathcal{K}^{(p)}_{0}(O_{X \times A} \ \mathrm{on} \ Y)$ be the eigenspace of $\psi^{m}=m^{p}$. It follows from Lemma \ref{Lemma: GilletSoule} that $L^{i,A}_{\bullet} \in \Gamma(U_{i}, \mathcal{K}^{(p)}_{0}(O_{X \times A} \ \mathrm{on} \ Y))$. Hence, $\alpha_{A}$ has the form
\[
\alpha_{A}: \mathrm{Hilb}(A) \to H^{0}(X, \mathcal{K}^{(p)}_{0}(O_{X \times A} \ \mathrm{on} \ Y)).
\]

According to Corti\~nas, Haesemeyer and Weibel \cite{CHW}, Chern character from K-theory to negative cyclic homology respects Adams operations. The natural map from negative cyclic homology to Hochschild homology also respects Adams operations, so does the map $ch_{A}$ (3.4)
\begin{equation*}
ch_{A}: H^{0}(X, \mathcal{K}^{(p)}_{0}(O_{X \times A} \ \mathrm{on} \ Y)) \to H^{0}(X, \mathcal{HH}^{(p)}_{0}(O_{X \times A} \ \mathrm{on} \ Y)).
\end{equation*}
The natural projection $r_{A}$ (3.5) also respects Adams operations
{\small 
 \begin{equation*}
  r_{A}: H^{0}(X, \mathcal{HH}^{(p)}_{0}(O_{X \times A} \ \mathrm{on} \  Y)) \to  H^{0}(X, \mathcal{\overline{HH}}^{(p)}_{0}(O_{X \times A} \ \mathrm{on} \  Y)),
 \end{equation*}
 }where $\mathcal{\overline{HH}}^{(p)}_{0}(O_{X \times A} \ \mathrm{on} \  Y)$ is the kernel of the map (induced by augmentation $A \to k$)
 \[
 \mathcal{HH}^{(p)}_{0}(O_{X \times A} \ \mathrm{on} \  Y) \to \mathcal{HH}^{(p)}_{0}(O_{X} \ \mathrm{on} \  Y).
 \] 
To summarize,
\begin{lemma} \label{Lemma: Eigenspace-Hilb-Map-HH}
With notation as above, the set-theoretic map $r_{A} \circ ch_{A} \circ \alpha_{A}: \mathrm{\mathbb{H}ilb}(A) \to H^{0}(X, \mathcal{\overline{HH}}_{0}(O_{X \times A} \ \mathrm{on} \ Y))$ defined in Definition \ref{Definition: Hilb-Map-Relative HH} has the form 
{\footnotesize
\begin{align*}
 & \mathrm{\mathbb{H}ilb}(A) \xrightarrow{\alpha_{A}} H^{0}(X,  \mathcal{K}^{(p)}_{0}(O_{X \times A} \ \mathrm{on} \  Y)) \xrightarrow{ch_{A}} H^{0}(X, \mathcal{HH}^{(p)}_{0}(O_{X \times A} \ \mathrm{on} \ Y)) \\ 
 & \xrightarrow{r_{A}} H^{0}(X, \mathcal{\overline{HH}}^{(p)}_{0}(O_{X \times A} \ \mathrm{on} \  Y)).
\end{align*}
}
\end{lemma}

We use the group $H^{0}(X, \mathcal{\overline{HH}}^{(p)}_{0}(O_{X \times A} \ \mathrm{on} \  Y))$ to define a functor.
\begin{definition} \label{Definition: Functor-HH-final}
One defines a functor
\begin{align*}
\mathbb{HH}^{(p)}: \ & Art_{k} \longrightarrow  Set_{*} \\
 & \ A \ \  \longrightarrow H^{0}(X, \mathcal{\overline{HH}}^{(p)}_{0}(O_{X \times A} \ \mathrm{on} \  Y)),
\end{align*}
where one forgets the group structure of $H^{0}(X, \mathcal{\overline{HH}}^{(p)}_{0}(O_{X \times A} \ \mathrm{on} \  Y))$ and considers it as a set (with zero as the chosen point).
\end{definition}

\begin{lemma} \label{Lemma: HH-final is Artin}
The functor $\mathbb{HH}^{(p)}$ is a functor of Artin rings. 
\end{lemma}
\begin{proof}
It is obvious that $\mathbb{HH}^{(p)}(k)=0$. The Hochschild homology is covariant, so is the functor $\mathbb{HH}^{(p)}$.
\end{proof}

By Theorem \ref{Theorem: Compute-HH-2-p case}, the morphism {\footnotesize$\mathcal{HH}^{(p)}_{0}(O_{X \times A} \ \mathrm{on} \  Y) \to \mathcal{HH}^{(p)}_{0}(O_{X} \ \mathrm{on} \  Y)$} has the form 
\[
\bigoplus_{j_{1}+j_{2}=p}\mathcal{H}^{p}_{Y}(\Omega^{j_{1}}_{X/k}) \otimes_{k} \Omega^{j_{2}}_{A/k}  \to \mathcal{H}^{p}_{Y}(\Omega^{p}_{X/k}) \otimes_{k} k.
\]

When $j_{2}=0$, $\Omega^{j_{2}}_{A/k}=\Omega^{0}_{A/k}=A=k \oplus m_{A}$ (as $k$-vector spaces), where $m_{A}$ is the maximal ideal of $A$. It follows that
{\small
\[
\mathcal{H}^{p}_{Y}(\Omega^{p}_{X/k}) \otimes_{k} \Omega^{0}_{A/k} =(\mathcal{H}^{p}_{Y}(\Omega^{p}_{X/k}) \otimes_{k} k) \bigoplus (\mathcal{H}^{p}_{Y}(\Omega^{p}_{X/k}) \otimes_{k} m_{A}).
\]
}
This implies that 
\begin{lemma} \label{Lemma: HH-is localCoh}
With notation as above, $\mathcal{\overline{HH}}^{(p)}_{0}(O_{X \times A} \ \mathrm{on} \  Y)$ is 
\[
 (\mathcal{H}^{p}_{Y}(\Omega^{p}_{X/k}) \otimes_{k} m_{A}) \bigoplus (\bigoplus_{ \tiny \begin{array}{c}
j_{1}+j_{2}=p \\
j_{2}>0
\end{array}} \mathcal{H}^{p}_{Y}(\Omega^{j_{1}}_{X/k}) \otimes_{k} \Omega^{j_{2}}_{A/k}).
 \]
Consequently, $\mathbb{HH}^{(p)}(A)$ is
 {\small
 \[
 (H^{0}(X, \mathcal{H}^{p}_{Y}(\Omega^{p}_{X/k})) \otimes_{k} m_{A}) \bigoplus (\bigoplus_{ \tiny \begin{array}{c}
j_{1}+j_{2}=p \\
j_{2}>0
\end{array}} H^{0}(X, \mathcal{H}^{p}_{Y}(\Omega^{j_{1}}_{X/k})) \otimes_{k} \Omega^{j_{2}}_{A/k}).
 \]
 }
\end{lemma}

\begin{remark} \label{Remark:HH-with-omegaXA}
Let $\overline{\Omega}^{p}_{X \times A/k}$ be the kernel of the morphism $\Omega^{p}_{X \times A/k} \to \Omega^{p}_{X/k}$ induced by the augmentation map $A \to k$. Using the isomorphism (3.10), we deduce that 
\[
\overline{\Omega}^{p}_{X \times A/k}= (\Omega^{p}_{X/k} \otimes_{k} m_{A}) \bigoplus (\bigoplus_{ \tiny \begin{array}{c}
j_{1}+j_{2}=p \\
j_{2}>0
\end{array}} \Omega^{j_{1}}_{X/k} \otimes_{k} \Omega^{j_{2}}_{A/k}).
\] 
By Lemma \ref{Lemma: HH-is localCoh}, $\mathcal{\overline{HH}}^{(p)}_{0}(O_{X \times A} \ \mathrm{on} \  Y)$ is isomorphic to the local cohomolgy sheaf $\mathcal{H}^{p}_{Y}(\overline{\Omega}^{p}_{X \times A/k})$ and there is an isomorphism 
\[
\mathbb{HH}^{(p)}(A) = H^{0}(X, \mathcal{H}^{p}_{Y}(\overline{\Omega}^{p}_{X \times A/k})).
\]
\end{remark}

Let $f: B \to A$ be a morphism in the category $Art_{k}$, there exists a commutative diagram (following from commutative diagrams (3.6) and (3.7) on page 14)
\[
\begin{CD}
 \mathrm{\mathbb{H}ilb}(B) @> r_{B} \circ ch_{B} \circ \alpha_{B}>>  \mathbb{HH}^{(p)}(B) \\
 @Vf_{Hilb}VV  @V\overline{f}_{HH}VV \\
  \mathrm{\mathbb{H}ilb}(A) @>r_{A} \circ ch_{A} \circ \alpha_{A}>> \mathbb{HH}^{(p)}(A).
 \end{CD}
\]

This shows that
\begin{lemma} \label{Lemma: transf-Hilb-HH-final}
With notation as above, there exists a natural transformation between functors of Artin rings
\begin{equation*}
\mathbb{T}:  \mathrm{\mathbb{H}ilb} \to \mathbb{HH}^{(p)},
\end{equation*}
which is defined to be, for $ A \in Art_{k}$, 
\[
\mathbb{T}(A)=r_{A} \circ ch_{A} \circ \alpha_{A}: \mathrm{\mathbb{H}ilb}(A) \to \mathbb{HH}^{(p)}(A).
\]
\end{lemma}

A natural transformation between two functors of Artin rings induces a map of the universal obstruction spaces (see the Remark on page 550 of \cite{FM-Obs}). Hence, the above natural transformation $\mathbb{T}$ induces a map of the universal obstruction spaces (as pointed sets)
\begin{equation}
O_{\mathrm{\mathbb{H}ilb}} \to O_{\mathbb{HH}^{(p)}}.
\end{equation}

\begin{lemma} \label{Lemma: HH-surjective}
Let $e: 0\to M \to B \xrightarrow{f} A \to 0$ be any small extension, the map $\mathbb{HH}^{(p)}(B) \to \mathbb{HH}^{(p)}(A)$ induced by $f$ is surjective.
\end{lemma}

\begin{proof}
Since $f: B \rightarrow A$ is surjective, for each integer $i \geq 0$, the induced map $\Omega^{i}_{B/k} \rightarrow \Omega^{i}_{A/k}$ is a  surjective morphism of $k$-vector spaces.  Since $f$ is a local morphism, the induced map $m_{B} \rightarrow m_{A}$ is surjective.

It follows from Lemma \ref{Lemma: HH-is localCoh} that $\mathbb{HH}^{(p)}(B) \to \mathbb{HH}^{(p)}(A)$ is surjective.

\end{proof}

This implies that the functor $\mathbb{HH}^{(p)}$ is smooth in the sense of Definition \ref{Definition: smooth functor}.
Consequently, by Lemma \ref{Lemma: smooth functor&trivial obs}, the universal obstruction theory of the functor $\mathbb{HH}^{(p)}$ is trivial, which is complete and linear. 
\begin{corollary} \label{Corollary:HH-obs space-trivial}
The universal obstruction space $O_{\mathbb{HH}^{(p)}}$ is a trivial space.
\end{corollary}

By Corollary \ref{Corollary: Obs HIlb is VS}, $O_{\mathrm{\mathbb{H}ilb}}$ is a vector space. Hence, the map (3.11) is a trivial morphism of vector spaces. In summary,
\begin{theorem} \label{Theorem: T induces trivial map}
The natural transformation $\mathbb{T}$ induces a trivial morphism of the universal obstruction spaces
\begin{equation}
O_{\mathrm{\mathbb{H}ilb}} \to O_{\mathbb{HH}^{(p)}},
\end{equation}
where both $O_{\mathrm{\mathbb{H}ilb}}$ and $O_{\mathbb{HH}^{(p)}}$ are vector spaces.
\end{theorem}

By Remark \ref{Remark: HH(A)-is-omegaXA} and Remark \ref{Remark:HH-with-omegaXA}, the map $r_{A} \circ ch_{A} \circ \alpha_{A}$ in Lemma \ref{Lemma: Eigenspace-Hilb-Map-HH} can be written as 
{\small
\begin{align}
 & \mathrm{\mathbb{H}ilb}(A) \xrightarrow{\alpha_{A}} H^{0}(X,  \mathcal{K}^{(p)}_{0}(O_{X \times A} \ \mathrm{on} \  Y)) \xrightarrow{ch_{A}} H^{0}(X, \mathcal{H}_{Y}^{p}(\Omega^{p}_{X \times A/ k}))  \\ 
 & \xrightarrow{r_{A}} H^{0}(X, \mathcal{H}_{Y}^{p}(\overline{\Omega}^{p}_{X \times A/ k})). \notag
\end{align}
} For later purpose, we need to describe the map $r_{A} \circ ch_{A} \circ \alpha_{A}$.

We first use a construction of Ang\'eniol and Lejeune-Jalabert \cite{A-LJ} to describe the map $ch_{A}$\footnote{When $A=k[\varepsilon]$, we have used this construction to describe a map from the tangent space $\mathrm{T_{Y}Hilb}^{p}(X)$ of the Hilbert scheme at the point $Y$ to local cohomology in section 3 of \cite{Y-3}.}. An element of $H^{0}(X,  \mathcal{K}^{(p)}_{0}(O_{X \times A} \ \mathrm{on} \  Y))$ is a cocycle $(C_{i})_{i}$, where $C_{i} \in \Gamma(U_{i},  \mathcal{K}^{(p)}_{0}(O_{X \times A} \ \mathrm{on} \  Y))$ is represented by a strict perfect complex $F^{i}_{\bullet}$
\[
 \begin{CD}
  0 @>>> L_{n} @>M_{n}>> L_{n-1} @>M_{n-1}>>  \cdots @>M_{2}>> L_{1} @>M_{1}>> L_{0} @>>> 0,
 \end{CD}
\]
whose homology is supported on $Y \cap U_{i}$. Here each $L_{l}$ is a free $O_{X \times A}(U_{i})$-modules of finite rank and each $M_{l}$ is a matrix with entries in $O_{X \times A}(U_{i})$.

\begin{definition} [Page 24 in \cite{A-LJ}] \label{definition: local cycles}
Let $p$ be any positive integer, the local fundamental class attached to this perfect complex $F^{i}_{\bullet}$ is defined to be the following collection
\[
 [F^{i}_{\bullet}]_{loc}=\{\frac{1}{p!}dM_{l}\circ dM_{l+1}\circ \dots \circ dM_{l+p-1}\}, l =  1, 2, \cdots,
\]
where $d=d_{k}$ and each $dM_{l}$ is the matrix of differentials. In other words, $dM_{l} \in \mathrm{Hom}(L_{l},L_{l-1}\otimes \Omega_{O_{X \times A}(U_{i}) /k}^{1})$.
\end{definition}

By Lemma 3.1.1 (on page 24) and Definition 3.4 (on page 29) in \cite{A-LJ}, the local fundamental class $[F^{i}_{\bullet}]_{loc}$ defines a cocycle of the complex $\mathcal{H}om(F^{i}_{\bullet}, \Omega^{p}_{O_{X \times A}(U_{i}) /k}\otimes F^{i}_{\bullet})$ and its image (still denoted $[F^{i}_{\bullet}]_{loc}$) in $\mathcal{E}XT^{p}(F^{i}_{\bullet}, \Omega^{p}_{O_{X \times A}(U_{i}) /k}\otimes F^{i}_{\bullet})$, which is the $p$-th cohomology of the complex $\mathcal{H}om(F^{i}_{\bullet}, \Omega^{p}_{O_{X \times A}(U_{i}) /k}\otimes F^{i}_{\bullet})$, does not depend on the choice of the basis of $F^{i}_{\bullet}$.

Since $F^{i}_{\bullet}$ is supported on $Y \cap U_{i}$, by the discussion after Definition 2.3.1 on page 98-99 in [1], there exists a trace map
\[
 \mathrm{Tr}:  \mathcal{E}XT^{p}(F^{i}_{\bullet}, \Omega^{p}_{O_{X \times A}(U_{i}) /k}\otimes F^{i}_{\bullet}) \to \Gamma(U_{i}, \mathcal{H}^{p}_{Y}(\Omega^{p}_{X \times A /k})).
\]

\begin{definition}   [Definition 2.3.2 on page 99 in \cite{A-LJ}] \label{Definition: Newton class}
The image of $[F^{i}_{\bullet}]_{loc}$ under the above trace map $\mathrm{Tr}$, denoted $\mathcal{V}_{F^{i}_{\bullet}}$, is called Newton class.
\end{definition}

The Grothendieck group of a triangulated category is the monoid of isomorphism objects modulo the submonoid formed from distinguished triangles.
\begin{lemma} [Proposition 4.3.1 on page 113 in \cite{A-LJ}] \label{Lemma: Newton class-Well}
The Newton class $\mathcal{V}_{F^{i}_{\bullet}}$ is well-defined on the Grothendieck group $K_{0}(O_{X \times A}(U_{i}) \ \mathrm{on} \ Y \cap U_{i})$.
\end{lemma}

\begin{definition} \label{Definition: Chern-Newton}
One uses Newton class $\mathcal{V}_{F^{i}_{\bullet}}$ to defines a morphism
\begin{align*}
\rho_{i}:\Gamma(U_{i},  \mathcal{K}^{(p)}_{0}&(O_{X \times A} \ \mathrm{on} \  Y))  \to  \Gamma(U_{i}, \mathcal{H}^{p}_{Y}(\Omega^{p}_{X \times A /k}))  \\
& \ \ \ F^{i}_{\bullet} \ \ \ \ \ \ \ \ \ \  \longrightarrow  \ \ \ \ \ \ \  \mathcal{V}_{F^{i}_{\bullet}}.
\end{align*}

\end{definition}
These morphisms $\rho_{i}$'s patch to give a morphism, which describes the map $ch_{A}$ in (3.13).

\begin{remark} \label{Remark: EXt-limit-Local}
Let $X$ be a noetherian scheme and let $Y \subset X$ be a closed subscheme with ideal sheaf $J$, for $\mathcal{F}$ a quasi-coherent $O_{X}$-module, there is an isomorphism
\[
H_{Y}^{p}(X, \mathcal{F}) = \varinjlim_{n \to \infty}Ext^{p}(O_{X}/J^{n}, \mathcal{F}).
\]
It follows that, for $\beta \in Ext^{p}(O_{X}/J, \mathcal{F})$, $\beta$ defines an element $[\beta] \in H_{Y}^{p}(X, \mathcal{F})$. We call $[\beta]$ the limit of $\beta$.

\end{remark}

Now, we are ready to give an description of the composition (3.13). In notation of Setting \ref{Setting:S}, for any element $Y' \in \mathrm{\mathbb{H}ilb}(A)$, $Y' \cap U_{i}$ is given by a regular sequence $f^{A}_{i1}, \cdots, f^{A}_{ip}$ of $O_{X \times A}(U_{i})$. Let $L^{i,A}_{\bullet}$ be the Koszul complex of the regular sequence $f^{A}_{i1}, \cdots, f^{A}_{ip}$. By Definition \ref{Definition:Hilb-Map-K}, $\alpha_{A}(Y')=(L^{i,A}_{\bullet})_{i}$.

Let $L^{i}_{\bullet}$ denote the Koszul complex associated to the regular sequence $f_{i1}, \cdots, f_{ip}$, which is a resolution of $O_{X}(U_{i})/(f_{i1}, \cdots, f_{ip})$ and has the form
\[
 L^{i}_{\bullet}: \ 0 \to L^{i}_{p} \to \cdots \to L^{i}_{0} \to 0,
\]
where each $L^{i}_{l}$ is defined as usual. In particular, $L^{i}_{p}=\wedge^{p}(O_{X}(U_{i})^{\oplus p}) \cong O_{X}(U_{i})$ and $L^{i}_{0}=O_{X}(U_{i})$.

Let $\omega^{A}_{i}=df^{A}_{i1}  \wedge \cdots \wedge df^{A}_{ip}$, the following diagram (denoted $\beta^{A}_{i}$)
{\tiny
\[
\begin{cases}
 \begin{CD}
    L^{i}_{\bullet} @>>>   O_{X}(U_{i})/(f_{i1}, \cdots, f_{ip}) @>>> 0\\
 L^{i}_{p}(\cong O_{X}(U_{i})) @>\omega^{A}_{i}>> L^{i}_{0} \otimes \Omega^{p}_{X\times A/k}(U_{i})(\cong \Omega^{p}_{X\times A/k}(U_{i})),
 \end{CD}
 \end{cases}
\]}defines an element in $Ext^{p}(O_{X}(U_{i})/(f_{i1}, \cdots, f_{ip}), \Omega^{p}_{X \times A/k}(U_{i}))$. The limit $[\beta^{A}_{i}] \in \Gamma(U_{i}, \mathcal{H}_{Y}^{p}(\Omega^{p}_{X \times A/k}))$ (see Remark \ref{Remark: EXt-limit-Local}) of $\beta^{A}_{i}$ is the Newton class $\mathcal{V}_{L^{i,A}_{\bullet}}$ attached to the Koszul complex $L^{i,A}_{\bullet}$, i.e., $ch_{A}(L^{i,A}_{\bullet})=[\beta^{A}_{i}]$.

Let $\overline{\omega}^{A}_{i}$ be the image of $\omega^{A}_{i}$ under the morphism $\Gamma(U_{i}, \Omega^{p}_{X \times A/k}) \to \Gamma(U_{i}, \overline{\Omega}^{p}_{X \times A/k})$, where $\overline{\Omega}^{p}_{X \times A/k}$ is the kernel of $\Omega^{p}_{X \times A/k} \to \Omega^{p}_{X/k}$ as in Remark \ref{Remark:HH-with-omegaXA}. Concretely, $\overline{\omega}^{A}_{i}= \omega^{A}_{i}-df_{i1}  \wedge \cdots \wedge df_{ip}$. The following diagram (denoted $\overline{\beta}^{A}_{i}$)
{\tiny
\[
\begin{cases}
 \begin{CD}
    L^{i}_{\bullet} @>>>   O_{X}(U_{i})/(f_{i1}, \cdots, f_{ip}) @>>> 0\\
  L^{i}_{p}(\cong O_{X}(U_{i})) @>\overline{\omega}^{A}_{i}>> L^{i}_{0} \otimes \overline{\Omega}^{p}_{X\times A/k}(U_{i})(\cong \overline{\Omega}^{p}_{X\times A/k}(U_{i})),
 \end{CD}
 \end{cases}
\]}defines an element in $Ext^{p}(O_{X}(U_{i})/(f_{i1}, \cdots, f_{ip}), \overline{\Omega}^{p}_{X \times A/k}(U_{i}))$.  The limit $[\overline{\beta}^{A}_{i}] \in \Gamma(U_{i}, \mathcal{H}_{Y}^{p}(\overline{\Omega}^{p}_{X \times A/k}))$ of $\overline{\beta}^{A}_{i}$  is $r_{A}([\beta^{A}_{i}])$, i.e.,  $r_{A} \circ ch_{A}(L^{i,A}_{\bullet})=[\overline{\beta}^{A}_{i}]$.

\begin{lemma} \label{Lemma: rcha(Y')}
With notation as above, for any $Y' \in \mathrm{\mathbb{H}ilb}(A)$, $r_{A} \circ ch_{A} \circ \alpha_{A}(Y')$ can be given by a cocycle $([\overline{\beta}^{A}_{i}])_{i}$:
\[
r_{A} \circ ch_{A} \circ \alpha_{A}(Y') = ([\overline{\beta}^{A}_{i}])_{i}.
\]
\end{lemma}

\section{Semi-regularity map and obstructions}

In this section, we reconstruct the semi-regularity map and prove Conjecture \ref{Conjecture:main}. We keep the notation of Setting \ref{Setting:S} and also write $U_{i} \times k[\varepsilon]=U_{i} \times _{\mathrm{Spec}(k)} \mathrm{Spec}(k[\varepsilon])$.

\subsection{Reconstruct semi-regularity map}
The following Lemma is used below. 
\begin{lemma} [page 59 of \cite{Bloch1}] \label{lemma: Bloch-local-Coh-Equal}
For each integer $m \geq 0$, there is an isomorphism
\begin{equation}
H^{m}(X, \mathcal{H}^{p}_{Y}(\Omega^{p-1}_{X/k}))=H^{p+m}_{Y}(X, \Omega^{p-1}_{X/k}).
\end{equation}
where $\mathcal{H}^{p}_{Y}(\Omega^{p-1}_{X/k})$ is the local cohomology sheaf.
\end{lemma}

Both the infinitesimal Abel-Jacobi map and the semi-regularity map can be described in terms of local cohomology. Let $\mathcal{I}$ be the ideal sheaf of $Y$ in $X$, the normal bundle $N_{Y/X}$ is defined to be $\mathcal{H}om_{Y}(\mathcal{I}/\mathcal{I}^{2}, O_{Y})$. According to Bloch \cite{Bloch1} (page 61), there is a morphism of sheaves 
\begin{equation}
\phi: N_{Y/X} \to \mathcal{H}^{p}_{Y}(\Omega^{p-1}_{X/k})
\end{equation}
defined by contraction with the cycle class. Concretely, given $\mu \in \Gamma(U_{i}, N_{Y/X})$, the morphism $\phi(\mu)$ can be described as (see equality (3) on page 201 of \cite{BF})
\begin{equation*}
\phi(\mu)= \sum^{l=p}_{l=1} (-1)^{l-1}\tilde{\mu}(f_{il})\dfrac{df_{i1} \wedge \cdots \wedge \hat{df_{il}} \wedge \cdots df_{ip} }{f_{i1}\cdots f_{ip}},
\end{equation*}
where $\tilde{\mu}(f_{il}) \in O_{X}(U_{i})$ is a lifting of $\mu(f_{il}) \in O_{Y}(U_{i})$ and $\hat{df_{il}}$ means to omit $df_{il}$. It can be checked that $\phi(\mu)$ gives a well-defined class in $\Gamma(U_{i}, \mathcal{H}^{p}_{Y}(\Omega^{p-1}_{X/k}))$.

For each integer $m \geq 0$, the morphism $\phi$ induces maps (still denoted $\phi$) on cohomology groups
\[
H^{m}(N_{Y/X}) \xrightarrow{\phi} H^{m}(X, \mathcal{H}^{p}_{Y}(\Omega^{p-1}_{X/k})).
\]
Composing $\phi$ with the natural map $H^{p+m}_{Y}(X, \Omega^{p-1}_{X/k}) \xrightarrow{L} H^{p+m}(X, \Omega^{p-1}_{X/k})$, one obtains the composition\footnote{Here we use the isomorphism (4.1).} 
\[
H^{m}(N_{Y/X}) \xrightarrow{\phi} H^{m}(X, \mathcal{H}^{p}_{Y}(\Omega^{p-1}_{X/k})) \xrightarrow{L} H^{p+m}(X, \Omega^{p-1}_{X/k}).
\]
When $m=0$ and $1$, the above composition gives the infinitesimal Abel-Jacobi map and the semi-regularity map respectively.

We reconstruct the semi-regularity map by reconstructing the morphism $\phi$ (4.2). The first step is to identify $\Gamma(U_{i}, N_{Y/X})$ with first order infinitesimal embedded deformations of $Y \cap U_{i}$ in $U_{i}\times k[\varepsilon]$.
\begin{lemma} [Theorem 2.4 in \cite{Hartshorne2}] \label{Lemma: Normal bundle-deform-correspond}
 The elements of $\Gamma(U_{i}, N_{Y/X})$ are naturally one-to-one correspondence with first order infinitesimal embedded deformations of  $Y \cap U_{i}$ in $U_{i}\times k[\varepsilon]$.
\end{lemma}

We sketch the one-to-one correspondence in Lemma \ref{Lemma: Normal bundle-deform-correspond} briefly. An element of $\Gamma(U_{i}, N_{Y/X})$ is given by a morphism $\mu: \mathcal{I}(U_{i}) \to O_{Y}(U_{i})$, which is determined by its images on $f_{i1}, \cdots, f_{ip}$. Let $\tilde{\mu}(f_{il}) \in O_{X}(U_{i})$ be a lifting of $\mu(f_{il})$, where $l=1, \cdots, p$, then the regular sequence $f_{i1}+ \varepsilon \tilde{\mu}(f_{i1}), \cdots, f_{ip}+ \varepsilon \tilde{\mu}(f_{ip})$ gives a first order infinitesimal embedded deformation $(Y \cap U_{i})^{'}$.

If $\tilde{\tilde{\mu}}(f_{il}) \in O_{X}(U_{i})$ is a different lifting of $\mu(f_{il})$, the regular sequence $f_{i1}+ \varepsilon \tilde{\tilde{\mu}}(f_{i1}) \cdots, f_{ip}+ \varepsilon \tilde{\tilde{\mu}}(f_{ip})$ gives the same infinitesimal embedded deformation $(Y \cap U_{i})^{'}$. Hence, there is a well-defined map 
\begin{align*}
\psi: \Gamma(U_{i}, \  N_{Y/X}&) \longrightarrow \mathrm{\mathbb{H}ilb}^{U_{i}}_{Y \cap U_{i}}(k[\varepsilon]) \\
& \mu \longrightarrow (Y \cap U_{i})^{'},
\end{align*}
where $\mathrm{\mathbb{H}ilb}^{U_{i}}_{Y \cap U_{i}}(k[\varepsilon])$ is the set of the first order infinitesimal
 embedded deformation of $Y \cap U_{i}$ in $U_{i} \times k[\varepsilon]$.

Let $X=U_{i}$ and $A=k[\varepsilon]$ in (3.13), the local version of the composition $r_{k[\varepsilon]} \circ ch_{k[\varepsilon]} \circ \alpha_{k[\varepsilon]}$ has the form
{\footnotesize
\begin{align*}
& \mathrm{\mathbb{H}ilb}^{U_{i}}_{U_{i}\cap Y}(k[\varepsilon]) \xrightarrow{\alpha_{k[\varepsilon]}} \Gamma(U_{i},  \mathcal{K}^{(p)}_{0}(O_{X \times k[\varepsilon]} \ \mathrm{on} \  Y))  \xrightarrow{ch_{k[\varepsilon]}}  \Gamma(U_{i},  \mathcal{H}^{p}_{Y}(\Omega^{p}_{X \times k[\varepsilon]/k})) \\
& \xrightarrow{r_{k[\varepsilon]}} \Gamma(U_{i},  \mathcal{H}^{p}_{Y}(\overline{\Omega}^{p}_{X \times k[\varepsilon]/k})),
\end{align*}
}
where $\overline{\Omega}^{p}_{X \times k[\varepsilon]/k}$ is the direct sum
\begin{equation}
 (\Omega^{p}_{X/k} \otimes_{k} (\varepsilon)) \oplus (\Omega^{p-1}_{X/k} \otimes_{k}\Omega^{1}_{k[\varepsilon]/k} )
\end{equation}
and $(\varepsilon)$ is the maximal ideal of $k[\varepsilon]$. Composing the map $\psi$ with $r_{k[\varepsilon]} \circ ch_{k[\varepsilon]} \circ \alpha_{k[\varepsilon]}$, we obtain a set-theoretic map
\begin{equation*}
\pi_{i}: \Gamma(U_{i}, \ N_{Y/X}) \to  \Gamma(U_{i},  \mathcal{H}^{p}_{Y}(\overline{\Omega}^{p}_{X \times k[\varepsilon]/k})).
\end{equation*}

\begin{lemma} \label{Lemma: Pi-group-map}
The map $\pi_{i}$ is a homomorphism of abelian groups.
\end{lemma}

\begin{proof}
Given an element $\mu \in \Gamma(U_{i}, N_{Y/X})$, $\psi(\mu)$ is a deformation of $Y \cap U_{i}$ given by the regular sequence $f_{i1}+\varepsilon \tilde{\mu}(f_{i1}), \cdots, f_{ip}+\varepsilon \tilde{\mu}(f_{ip})$, where each $\tilde{\mu}(f_{il}) \in O_{X}(U_{i})$ is a lifting of $\mu(f_{il})$. In notation of Lemma \ref{Lemma: rcha(Y')}, $\pi_{i}(\mu)=r_{k[\varepsilon]} \circ ch_{k[\varepsilon]} \circ \alpha_{k[\varepsilon]}( \psi(\mu))$ is the limit of the following diagram
\[
\begin{cases}
 \begin{CD}
   L^{i}_{\bullet} @>>> O_{X}(U_{i})/(f_{i1}\cdots f_{ip}) @>>> 0\\
  L^{i}_{p} @>\omega_{1} \varepsilon + \omega_{2} d \varepsilon>> L^{i}_{0} \otimes \overline{\Omega}^{p}_{X \times k[\varepsilon]/k}(U_{i}), 
 \end{CD}
 \end{cases}
\]
where $\omega_{1}$ and $\omega_{2}$ are as follows
 \begin{align*}
 \omega_{1} & =\sum^{l=p}_{l=1} df_{i1} \wedge \cdots \wedge df_{i(l-1)} \wedge d \tilde{\mu}(f_{il}) \wedge df_{i(l+1)} \wedge \cdots \wedge df_{ip}, \\
  \omega_{2}& =\sum^{l=p}_{l=1} (-1)^{p-l}\tilde{\mu}(f_{il})df_{i1} \wedge \cdots \wedge \hat{df_{il}} \wedge \cdots df_{ip}.
 \end{align*}

From this description of $\pi_{i}(\mu)$, we deduce that $\pi_{i}$ is a homomorphism of abelian groups.

\end{proof}

These morphisms $\pi_{i}$'s patch to give a morphism of sheaves,
\begin{equation*}
\pi: N_{Y/X}  \rightarrow \mathcal{H}^{p}_{Y}( \overline{\Omega}^{p}_{X \times k[\varepsilon]/k}).
\end{equation*}

The contraction map $\Omega_{k[\varepsilon]/k}^{1} \to k$ induces a morphism
\begin{equation*}
\rfloor: \mathcal{H}^{p}_{Y}(\overline{\Omega}^{p}_{X \times k[\varepsilon]/k}) \to  \mathcal{H}^{p}_{Y}(\Omega^{p-1}_{X/k}).
\end{equation*}
 
\begin{corollary} \label{Corollary: Pi-Equal-Bloch's Map}
The composition of morphisms of sheaves
\begin{equation*}
N_{Y/X} \xrightarrow{\pi}  \mathcal{H}^{p}_{Y}(\overline{\Omega}^{p}_{X \times k[\varepsilon]/k}) \xrightarrow{\rfloor}  \mathcal{H}^{p}_{Y}(\Omega^{p-1}_{X/k})
\end{equation*}
agrees with the morphism $N_{Y/X} \xrightarrow{\phi(4.2)} \mathcal{H}^{p}_{Y}(\Omega^{p-1}_{X})$
in the sense that
\[
\rfloor \circ \pi = (-1)^{p-1} \phi.
\]
In other words, when $p$ is odd, $\rfloor \circ \pi = \phi$; when $p$ is even, $\rfloor \circ \pi = (-1)\phi$.
\end{corollary}

\begin{proof}
In notation of the proof of Lemma \ref{Lemma: Pi-group-map}, the composition $\rfloor \circ \pi (\mu)$ locally on $U_{i}$ is given by the following diagram
\[
\begin{cases}
 \begin{CD}
   L^{i}_{\bullet} @>>> O_{X}(U_{i})/(f_{i1}\cdots f_{ip}) @>>> 0\\
  L^{i}_{p} @>\omega_{2}>> L^{i}_{0} \otimes \Omega_{X/k}^{p-1}(U_{i}).
 \end{CD}
 \end{cases}
\]
Comparing with the description of $\phi$ (on page 23), one sees that $\rfloor \circ \pi = (-1)^{p-1}\phi$.
\end{proof}

Since the infinitesimal Abel-Jacobi map is the composition
\[ 
H^{0}(Y, N_{Y/X}) \xrightarrow{\phi} H^{0}(X, \mathcal{H}^{p}_{Y}(\Omega^{p-1}_{X/k})) \xrightarrow{L} H^{p}(X, \Omega^{p-1}_{X/k}),
\]
it follows from Corollary \ref{Corollary: Pi-Equal-Bloch's Map} that the infinitesimal Abel-Jacobi map agrees with the composition
{ \footnotesize
\[
H^{0}(Y, N_{Y/X}) \xrightarrow{\pi} H^{0}(X, \mathcal{H}^{p}_{Y}(\overline{\Omega}^{p}_{X \times k[\varepsilon]/k})) \xrightarrow{\rfloor} H^{0}(X, \mathcal{H}^{p}_{Y}(\Omega^{p-1}_{X/k})) \xrightarrow{L} H^{p}(X, \Omega^{p-1}_{X/k}).
\]
}
\  It is well known that the Zariski tangent space $\mathrm{\mathbb{H}ilb}(k[\varepsilon])$ of the local Hilbert functor $\mathrm{\mathbb{H}ilb}$ can be identified with $H^{0}(Y,N_{Y/X})$. Moreover, it follows from Remark \ref{Remark:HH-with-omegaXA} that the tangent space $\mathrm{\mathbb{HH}^{(p)}}(k[\varepsilon])$ of the functor $\mathrm{\mathbb{HH}^{(p)}}$ is $H^{0}(X, \mathcal{H}^{p}_{Y}(\overline{\Omega}^{p}_{X \times k[\varepsilon]/k}))$, where $\overline{\Omega}^{p}_{X \times k[\varepsilon]/k}$ is (4.3). Hence, the morphism $\pi$ is a morphism between tangent spaces of functors of Artin rings.

To summarize, we have proved the following theorem which answers (after slight modification) part (I) of Question \ref{Question:2}.
\begin{theorem} \label{Theorem: Inf'l AJ-tang-space}
The infinitesimal Abel-Jacobi map agrees with the composition
{ \footnotesize
 \[
H^{0}(Y, N_{Y/X}) \xrightarrow{\pi} H^{0}(X, \mathcal{H}^{p}_{Y}(\overline{\Omega}^{p}_{X \times k[\varepsilon]/k})) \xrightarrow{\rfloor} H^{0}(X, \mathcal{H}^{p}_{Y}(\Omega^{p-1}_{X/k})) \xrightarrow{L} H^{p}(X, \Omega^{p-1}_{X/k}),
\]
}where $\pi$ is a map between tangent spaces of two functors of Artin rings $\mathrm{\mathbb{H}ilb}$ and $\mathrm{\mathbb{HH}^{(p)}}$.
\end{theorem}

Since the semi-regularity map is the composition
\[ 
H^{1}(Y, N_{Y/X}) \xrightarrow{\phi} H^{1}(X, \mathcal{H}^{p}_{Y}(\Omega^{p-1}_{X/k})) \xrightarrow{L} H^{p+1}(X, \Omega^{p-1}_{X/k}),
\]
we deduce from Corollary \ref{Corollary: Pi-Equal-Bloch's Map} that
\begin{lemma} \label{Lemma: Semi-regularity-Agree}
The semi-regularity map agrees with the composition
{\footnotesize
 \[
H^{1}(Y, N_{Y/X}) \xrightarrow{\pi} H^{1}(X, \mathcal{H}^{p}_{Y}(\overline{\Omega}^{p}_{X \times k[\varepsilon]/k})) \xrightarrow{\rfloor} H^{1}(X, \mathcal{H}^{p}_{Y}(\Omega^{p-1}_{X/k})) \xrightarrow{L} H^{p+1}(X, \Omega^{p-1}_{X/k}).
\]
}
\end{lemma}

\subsection{Obstruction space}
We give a different proof of Theorem \ref{t:Bloch-semi-conj} in this subsection. The notation of Setting \ref{Setting:S} is still used here. 

Let $e$ be a principal small extension as in Theorem \ref{t:Bloch-semi-conj}. That is, $e$ is a principal small extension 
\[
e: 0 \to (\eta) \to B \xrightarrow{f} A \to 0,
\]
where $\eta^2=0$ and the differential $d: (\eta) \to \Omega_{B/k} \otimes_{B}A$ is injective. By Definition \ref{Definition: Functor-HH-final}, an element $\alpha^{A} \in \mathbb{HH}^{(p)}(A)$ has the form $\alpha^{A}=(\alpha^{A}_{i})_{i}$, where $\alpha^{A}_{i}$ is an element of $\Gamma(U_{i}, \mathcal{\overline{HH}}^{(p)}_{0}(O_{X \times A} \ \mathrm{on} \  Y))$ satisfying that $\alpha^{A}_{i}=\alpha^{A}_{j}$ over the intersection $U_{ij}=U_{i} \cap U_{j}$.

By mimicking the proof of Lemma \ref{Lemma: HH-surjective}, one sees that the morphism $\Gamma(U_{i}, 
 \mathcal{\overline{HH}}^{(p)}_{0}(O_{X \times B} \ \mathrm{on} \  Y)) \to \Gamma(U_{i}, \mathcal{\overline{HH}}^{(p)}_{0}(O_{X \times A} \ \mathrm{on} \  Y))$ induced by $f$ is surjective. Hence, one can lift $\alpha^{A}_{i}$ to $\alpha^{B}_{i}\in \Gamma(U_{i}, \mathcal{\overline{HH}}^{(p)}_{0}(O_{X \times B} \ \mathrm{on} \  Y))$.

The kernel of the morphism (induced by $f$) $m_{B} \to m_{A}$ is $(\eta)$. For a positive integer $j_{2}$ satisfying $0<j_{2}\leq p$, let $W_{j_{2}}$ denote the kernel of the morphism of $k$-vector spaces $\Omega^{j_{2}}_{B/k} \to \Omega^{j_{2}}_{A/k}$ and let $\mathcal{W}$ be the sheaf
\begin{equation}
\mathcal{W}=(\Omega^{p}_{X/k} \otimes_{k} (\eta)) \bigoplus (\bigoplus_{ \tiny \begin{array}{c}
j_{1}+j_{2}=p \\
j_{2}>0
\end{array}} \Omega^{j_{1}}_{X/k} \otimes_{k} W_{j_{2}}).
\end{equation}
Let $\bar{f}: \Gamma(U_{ij}, \mathcal{\overline{HH}}^{(p)}_{0}(O_{X \times B} \ \mathrm{on} \  Y)) \to \Gamma(U_{ij}, \mathcal{\overline{HH}}^{(p)}_{0}(O_{X \times A} \ \mathrm{on} \  Y))$ be the morphism induced by $f$, it follows from Lemma \ref{Lemma: HH-is localCoh} that the kernel of $\bar{f}$ is $\Gamma(U_{ij}, \mathcal{H}^{p}_{Y}(\mathcal{W}))$, where $\mathcal{H}^{p}_{Y}(\mathcal{W})$ denotes the local cohomology sheaf.

The difference of two liftings $\alpha^{B}_{i}$ and $\alpha^{B}_{j}$ over the intersection $U_{ij}$ satisfies that
\[
\bar{f}(\alpha^{B}_{i} - \alpha^{B}_{j})=  \alpha^{A}_{i} - \alpha^{A}_{j} = 0.
\]
This shows that $\alpha^{B}_{i} - \alpha^{B}_{j}$ lies in the kernel of $\bar{f}$, i.e., $(\alpha^{B}_{i} - \alpha^{B}_{j}) \in \Gamma(U_{ij}, \mathcal{H}^{p}_{Y}(\mathcal{W}))$. On the intersection $U_{ijk}=U_{i} \bigcap U_{j} \bigcap U_{k}$,
\[
\alpha^{B}_{i} - \alpha^{B}_{k} = (\alpha^{B}_{i} - \alpha^{B}_{j}) + ( \alpha^{B}_{j} - \alpha^{B}_{k}),
\]
where $\alpha^{B}_{k} \in \Gamma(U_{k}, \mathcal{\overline{HH}}^{(p)}_{0}(O_{X \times B} \ \mathrm{on} \  Y))$ is a lifting of $\alpha^{A}_{k}$. This shows that $(\alpha^{B}_{i}- \alpha^{B}_{j})$ forms a $\check{\mathrm{C}}$ech 1-cocycle $(\alpha^{B}_{i}- \alpha^{B}_{j})_{i,j}$, which defines a cohomology class of $H^{1}(X, \mathcal{H}^{p}_{Y}(\mathcal{W}))$, denoted $\{\alpha^{B}_{i} - \alpha^{B}_{j}\}_{i,j}$.

If $\alpha^{'B}_{i} \in \Gamma(U_{i}, \mathcal{\overline{HH}}^{(p)}_{0}(O_{X \times B} \ \mathrm{on} \  Y))$ is a different lifting of $\alpha^{A}_{i}$, on the intersection $U_{i} \cap U_{j}$, 
\[
\alpha^{B}_{i} - \alpha^{B}_{j} = ( \alpha^{B}_{i} - \alpha^{'B}_{i}) + ( \alpha^{'B}_{i} - \alpha^{'B}_{j}) - ( \alpha^{B}_{j} - \alpha^{'B}_{j}).
\]
Since $(\alpha^{B}_{i} - \alpha^{'B}_{i}) \in \Gamma(U_{i}, \mathcal{H}^{p}_{Y}(\mathcal{W}))$ and $(\alpha^{B}_{j} - \alpha^{'B}_{j}) \in \Gamma(U_{j}, \mathcal{H}^{p}_{Y}(\mathcal{W}))$, $((\alpha^{B}_{i} - \alpha^{'B}_{i}) -  (\alpha^{B}_{j} - \alpha^{'B}_{j}))_{i,j}$ is a 1-coboundary. Hence, the $\check{\mathrm{C}}$ech 1-cocycles $(\alpha^{B}_{i} - \alpha^{B}_{j})_{i,j}$ and $(\alpha^{'B}_{i} - \alpha^{'B}_{j})_{i,j}$ define the same cohomology class of $H^{1}(X, \mathcal{H}^{p}_{Y}(\mathcal{W}))$. This shows that the cohomology class $\{\alpha^{B}_{i} - \alpha^{B}_{j}\}_{i,j}$ is independent of the choice of $\alpha^{B}_{i}$.

\begin{definition} \label{Definition: Delta-1}
With notation as above, one defines a map\footnote{We will show that this is a trivial map below.}
\begin{align*}
\delta_{1}: \mathbb{HH}^{(p)}(A) & \to H^{1}(X, \mathcal{H}^{p}_{Y}(\mathcal{W}))  \\
(\alpha^{A}_{i})_{i} & \to \{\alpha^{B}_{i} - \alpha^{B}_{j}\}_{i,j}.
\end{align*}
\end{definition}

It is obvious that $d\eta$ lies in the kernel of the morphism $ \Omega^{1}_{B/k} \to \Omega^{1}_{A/k}$, i.e., $d\eta \in W_{1}$. Since $d: (\eta) \to \Omega_{B/k} \otimes_{B}A$ is injective, $d\eta \neq 0$, then $d\eta$ generates a one-dimensional subspace of the $k$-vector space $W_{1}$, denoted $(d\eta)$.

Let $h: W_{1} \to (d\eta)$ be the projection of $k$-vector spaces. The morphisms  $(\eta) \xrightarrow{=} (\eta)$ and $\ W_{1} \xrightarrow{h} (d\eta)$ induce a morphism of sheaves
\begin{equation}
\mathcal{W} \to (\Omega^{p}_{X/k} \otimes_{k} (\eta))  \oplus (\Omega^{p-1}_{X/k} \otimes_{k} (d\eta)).
\end{equation}
Since both $(d\eta)$ and $(\eta)$ are one dimensional $k$-vector spaces, one identifies $(d\eta)$ with $(\eta)$ and sees that there is a morphism of sheaves
\[
\mathcal{W} \to V \otimes_{k} (\eta),
\]
where $V=\Omega^{p}_{X/k}  \bigoplus \Omega^{p-1}_{X/k}$, which induces a morphism on cohomology
\[
\delta_{2}: H^{1}(X, \mathcal{H}^{p}_{Y}(\mathcal{W})) \to H^{1}(X,  \mathcal{H}^{p}_{Y}(V)) \otimes_{k} (\eta).
\]

For $k[\varepsilon]$ the ring of dual numbers, both $\Omega^{1}_{k[\varepsilon]/k}$ and the maximal ideal $(\varepsilon)$ are one dimensional $k$-vector spaces. There is an isomorphism of sheaves $V \cong \overline{\Omega}^{p}_{X \times k[\varepsilon]/k}$, where $\overline{\Omega}^{p}_{X \times k[\varepsilon]/k}$ is (4.3), which induces an isomorphism on cohomology
\[
i: H^{1}(X,  \mathcal{H}^{p}_{Y}(V))\otimes_{k} (\eta) \to H^{1}(X, \mathcal{H}^{p}_{Y}(\overline{\Omega}^{p}_{X \times k[\varepsilon]/k}))\otimes_{k} (\eta).
\]

\begin{definition} \label{Definition: Delta}
With notation as above, we define a map 
\begin{equation}
\delta_{e}: \mathbb{HH}^{(p)}(A) \times \check{(\eta)}  \to H^{1}(X, \mathcal{H}^{p}_{Y}(\overline{\Omega}^{p}_{X \times k[\varepsilon]/k}))
\end{equation}
to be the composition
{\tiny
\[
 \mathrm{\mathbb{HH}}^{(p)}(A) \times \check{(\eta)}  \xrightarrow{(i \circ \delta_{2} \circ \delta_{1}, 1)} (H^{1}(X, \mathcal{H}^{p}_{Y}(\overline{\Omega}^{p}_{X \times k[\varepsilon]/k}))\otimes_{k} (\eta))\times \check{(\eta)} \xrightarrow{\tilde{\delta}} H^{1}(X, \mathcal{H}^{p}_{Y}(\overline{\Omega}^{p}_{X \times k[\varepsilon]/k})),
\]}where $\check{(\eta)}$ is the dual space of the $k$-vector space $(\eta)$ and $\tilde{\delta}(x \otimes \eta, g)=g(\eta)x$.

\end{definition}

\begin{lemma} \label{Lemma: Main Lemma-Commu-Diag}
There exists a commutative diagram of sets
\[
\begin{CD}
  \mathrm{\mathbb{H}ilb}(A)\times \check{(\eta)} @>(\mathbb{T}(A), 1)>> \mathbb{HH}^{(p)}(A) \times \check{(\eta)}\\
  @V v_{e} V(2.3)V  @V \delta_{e} V(4.6)V \\
H^{1}(Y, N_{Y/X}) @>\pi >>   H^{1}(X, \mathcal{H}^{p}_{Y}(\overline{\Omega}^{p}_{X \times k[\varepsilon]/k})), 
 \end{CD}
\]
where $\mathbb{T}(A)=r_{A} \circ ch_{A} \circ \alpha_{A}$ in Lemma \ref{Lemma: transf-Hilb-HH-final} and $\pi$ is the morphism in Lemma \ref{Lemma: Semi-regularity-Agree}.
\end{lemma}

\begin{proof}
For $Y' \in \mathrm{\mathbb{H}ilb}(A)$ and $g \in \check{(\eta)}$, we need to prove that $\pi \circ v_{e}(Y', g)= \delta_{e} \circ (\mathbb{T}(A),1)(Y', g)$.

In notation of Setting \ref{Setting:S}, the regular sequence $f^{A}_{i1}, \cdots, f^{A}_{ip}$ can be lifted to a regular sequence $f^{B}_{i1}, \cdots, f^{B}_{ip}$ of $O_{X}(U_{i}) \otimes_{k} B$. On the intersection $U_{ij}=U_{i} \cap U_{j}$, two liftings $f^{B}_{i1}, \cdots, f^{B}_{ip}$ and $f^{B}_{j1}, \cdots, f^{B}_{jp}$ satisfy that $f^{B}_{il}-f^{B}_{jl}=\eta h_{ijl}$, where $h_{ijl} \in O_{X}(U_{ij})$ and $l=1, \cdots,p$.

Let $h'_{ijl}$ be the image of $h_{ijl}$ in $O_{Y}(U_{ij})$. By the description of $v_{e}$ (2.4), $v_{e}(Y', g)$ is given by a $\check{\mathrm{C}}$ech 1-cocycle formed from $g(\eta)\mu_{ij}$, where $\mu_{ij} \in \Gamma(U_{ij},N_{Y/X})$ is the morphism mapping $f_{il}$ to $h'_{ijl}$, see page 8-9 for details.

Let $L^{ij}_{\bullet}$ be the Koszul  resolution of $O_{X}(U_{ij}) /(f_{i1}, \cdots, f_{ip})$, which has the form
\begin{equation*}
 L^{ij}_{\bullet}: \ 0 \to L^{ij}_{p} \to \cdots \to L^{ij}_{0} \to 0.
\end{equation*}
The description of $\pi \circ v_{e}(Y', g)$ has been essentially given in the proof of Lemma \ref{Lemma: Pi-group-map}, we recall it briefly. It is trivial that $h_{ijl}$ is a lifting of $\mu_{ij}(f_{il})=h'_{ijl}$ to $O_{X}(U_{ij})$, let $\omega^{1}_{ij}$ and $\omega^{2}_{ij} $ be as follows
 \begin{align}
 \omega^{1}_{ij} & =\sum^{l=p}_{l=1} df_{i1} \wedge \cdots \wedge df_{i(l-1)} \wedge dh_{ijl} \wedge df_{i(l+1)} \wedge \cdots \wedge df_{ip}, \\
  \omega^{2}_{ij}& =\sum^{l=p}_{l=1} (-1)^{p-l}h_{ijl}df_{i1} \wedge \cdots \wedge \hat{df_{il}} \wedge \cdots df_{ip}.
 \end{align} 
The following diagram (denoted $\overline{\beta}^{k[\varepsilon]}_{ij}$)
\[
\begin{cases}
 \begin{CD}
   L^{ij}_{\bullet} @>>> O_{X}(U_{ij}) /(f_{i1}, \cdots, f_{ip}) @>>> 0 \\
 L^{ij}_{p} @>g(\eta) (\omega^{1}_{ij}\varepsilon+\omega^{2}_{ij} d \varepsilon) >> L^{ij}_{0} \otimes 
 \overline{\Omega}^{p}_{X \times k[\varepsilon]/k}(U_{ij}),
 \end{CD}
 \end{cases}
\]
defines an element in $Ext^{p}(O_{X}(U_{ij})/(f_{i1}, \cdots, f_{ip}), \overline{\Omega}^{p}_{X \times k[\varepsilon]/k}(U_{ij}) )$. The $\mathrm{\check{C}}$ech 1-cocycles formed from the limit $[\overline{\beta}^{k[\varepsilon]}_{ij}] \in \Gamma(U_{ij}, \mathcal{H}_{Y}^{p}(\overline{\Omega}^{p}_{X \times k[\varepsilon]/k}))$ of $\overline{\beta}^{k[\varepsilon]}_{ij}$ (see Remark \ref{Remark: EXt-limit-Local}) gives a class of $H^{1}(X, \mathcal{H}^{p}_{Y}(\overline{\Omega}^{p}_{X \times k[\varepsilon]/k}))$, which is $\pi \circ v_{e}(Y', g)$.

In notation of  Lemma \ref{Lemma: rcha(Y')}, let $\overline{\omega}^{A}_{i}=df^{A}_{i1}  \wedge \cdots \wedge df^{A}_{ip}-df_{i1}  \wedge \cdots \wedge df_{ip}$, $\mathbb{T}(A)(Y')=([\overline{\beta}^{A}_{i}])_{i}$, where $[\overline{\beta}^{A}_{i}]$ is the limit of the following diagram
\[
\begin{cases}
 \begin{CD}
    L^{i}_{\bullet} @>>>   O_{X}(U_{i})/(f_{i1}, \cdots, f_{ip}) @>>> 0 \\
  L^{i}_{p} @>\overline{\omega}^{A}_{i}>> L^{i}_{0} \otimes \overline{\Omega}^{p}_{O_{X\times A}(U_{i}) / k}.
 \end{CD}
 \end{cases}
\]

Let $L^{i,B}_{\bullet}$ be the Koszul complex associated to the regular sequence $f^{B}_{i1}, \cdots, f^{B}_{ip}$. Let $B=A$ in Lemma \ref{Lemma: rcha(Y')}, the following diagram (denoted $\overline{\beta}^{B}_{i}$)
\[
\begin{cases}
 \begin{CD}
    L^{i}_{\bullet} @>>>   O_{X}(U_{i})/(f_{i1}, \cdots, f_{ip}) @>>> 0 \\
  L^{i}_{p} @>\overline{\omega}^{B}_{i}>> L^{i}_{0} \otimes \overline{\Omega}^{p}_{O_{X \times B}(U_{i}) / k},
 \end{CD}
 \end{cases}
\]
where $\overline{\omega}^{B}_{i}=df^{B}_{i1}  \wedge \cdots \wedge df^{B}_{ip}-df_{i1}  \wedge \cdots \wedge df_{ip}$, defines an element of $Ext^{p}(O_{X}(U_{i})/(f_{i1}, \cdots, f_{ip}), \overline{\Omega}^{p}_{X \times B/k}(U_{i}) )$. The limit $[\overline{\beta}^{B}_{i}]$ of $\overline{\beta}^{B}_{i}$ is $r_{B} \circ ch_{B}(L^{i,B}_{\bullet})$.

It is obvious that $[\overline{\beta}^{B}_{i}]$ lifts $[\overline{\beta}^{A}_{i}]$. To find a lifting of $\mathbb{T}(A)(Y')$, i.e., a lifting of $([\overline{\beta}^{A}_{i}])_{i}$, we need to check the difference $[\overline{\beta}^{B}_{i}]-[\overline{\beta}^{B}_{j}]$. It follows from Remark \ref{Remark:HH-with-omegaXA} that the kernel of the morphism $\overline{\Omega}^{p}_{X \times B/k} \to \overline{\Omega}^{p}_{X \times A/k}$ is the sheaf $\mathcal{W}$ (4.4). Over the intersection $U_{ij}$, $\overline{\omega}^{B}_{i}-\overline{\omega}^{B}_{j}$ lies in the kernel of the morphism $\Gamma(U_{ij}, \overline{\Omega}^{p}_{X \times B/k}) \to \Gamma(U_{ij}, \overline{\Omega}^{p}_{X \times A/k})$ induced by $B \xrightarrow{f} A$, which is $\Gamma(U_{ij}, \mathcal{W})$, i.e., $\overline{\omega}^{B}_{i}-\overline{\omega}^{B}_{j} \in \Gamma(U_{ij}, \mathcal{W})$. Hence, $[\overline{\beta}^{B}_{i}]-[\overline{\beta}^{B}_{j}]$ lies in $\Gamma(U_{ij}, \mathcal{H}_{Y}^{p}(\mathcal{W}))$. The $\check{\mathrm{C}}$ech 1-cocycle formed by $[\overline{\beta}^{B}_{i}]-[\overline{\beta}^{B}_{j}]$ defines a class of $H^{1}(X, \mathcal{H}_{Y}^{p}(\mathcal{W}))$, which is $\delta_{1} \circ \mathrm{T}(A)(Y')$.

The morphisms (4.5) induces a morphism
\[
Q: \Gamma(U_{ij}, \mathcal{W}) \to \Gamma(U_{ij}, (\Omega^{p}_{X /k} \otimes_{k} (\eta)) \oplus (\Omega^{p-1}_{X /k} \otimes_{k}(d\eta))).
\] We want to describe the image $Q(\overline{\omega}^{B}_{i}-\overline{\omega}^{B}_{j})$. Since $f^{B}_{il}-f^{B}_{jl}= \eta h_{ijl}$, where $l=1, \cdots, p$, it follows that $df^{B}_{il}-df^{B}_{jl} = \eta dh_{ijl} + h_{ijl} d\eta$. Since $\eta^{2}=0$, $\eta d\eta=0$. It follows that
{\tiny
\begin{align*}
\overline{\omega}^{B}_{i} & = df^{B}_{i1} \wedge \cdots \wedge df^{B}_{ip}- df_{i1} \wedge \cdots \wedge df_{ip} \\
&=(df^{B}_{j1} + \eta dh_{ij1} + h_{ij1} d\eta) \wedge \cdots \wedge (df^{B}_{jp} + \eta dh_{ijp} + h_{ijp} d\eta) - df_{i1} \wedge \cdots \wedge df_{ip} \\
& = df^{B}_{j1} \wedge \cdots \wedge df^{B}_{jp} + (\sum^{l=p}_{l=1} df^{B}_{i1} \wedge \cdots \wedge df^{B}_{i(l-1)} \wedge dh_{ijl} \wedge df^{B}_{i(l+1)} \wedge \cdots \wedge df^{B}_{ip}) \eta  \\
& + (\sum^{l=p}_{l=1} (-1)^{p-l}h_{ijl}df^{B}_{i1} \wedge \cdots \wedge \hat{df^{B}_{il}} \wedge \cdots df^{B}_{ip}) d\eta - df_{i1} \wedge \cdots \wedge df_{ip} \\
& = \overline{\omega}^{B}_{j} + (\sum^{l=p}_{l=1} df^{B}_{i1} \wedge \cdots \wedge df^{B}_{i(l-1)} \wedge dh_{ijl} \wedge df^{B}_{i(l+1)} \wedge \cdots \wedge df^{B}_{ip}) \eta  \\
& + (\sum^{l=p}_{l=1} (-1)^{p-l}h_{ijl}df^{B}_{i1} \wedge \cdots \wedge \hat{df^{B}_{il}} \wedge \cdots df^{B}_{ip}) d\eta.
 \end{align*}    
 } This implies that
\[
 Q(\overline{\omega}^{B}_{i}-\overline{\omega}^{B}_{j})= \omega^{1}_{ij}\eta + \omega^{2}_{ij} d \eta,
\]
where $\omega^{1}_{ij}$ is (4.7) and $\omega^{2}_{ij}$ is (4.8).

Let $V=\Omega^{p}_{X /k} \oplus \Omega^{p-1}_{X /k}$ as on page 28. The following diagram (denoted $\beta'_{ij}$)
\[
\begin{cases}
 \begin{CD}
   L^{ij}_{\bullet} @>>> O_{X}(U_{ij}) /(f_{i1}, \cdots, f_{ip}) @>>> 0\\
 L^{ij}_{p} @>(\omega^{1}_{ij}+\omega^{2}_{ij}) \eta>> L^{ij}_{0} \otimes V(U_{ij})\otimes(\eta),
 \end{CD}
 \end{cases}
\]
defines an element in $Ext^{p}(O_{X}(U_{ij})/(f_{i1}, \cdots, f_{ip}), V(U_{ij})\otimes(\eta))$. The limit $[\beta'_{ij}]$ of $\beta'_{ij}$ is in $\Gamma(U_{ij}, \mathcal{H}_{Y}^{p}(V)\otimes(\eta))$. The $\check{\mathrm{C}}$ech 1-cocycle formed by $[\beta'_{ij}]$ defines a class of $H^{1}(X, \mathcal{H}_{Y}^{p}(V))\otimes(\eta)$, which is $\delta_{2} \circ \delta_{1}(\mathbb{T}(A)(Y'))$.

The isomorphism $V \cong \overline{\Omega}^{p}_{X \times k[\varepsilon]/k}$ maps $\omega^{1}_{ij}+ \omega^{2}_{ij}$ to $\omega^{1}_{ij} \varepsilon + \omega^{2}_{ij} d \varepsilon$, so $i \circ \delta_{2} \circ \delta_{1}(\mathbb{T}(A)(Y'))$ is given by the $\check{\mathrm{C}}$ech 1-cocycle formed by the limit of the diagram
\[
\begin{cases}
 \begin{CD}
   L^{ij}_{\bullet} @>>> O_{X}(U_{ij}) /(f_{i1}, \cdots, f_{ip}) @>>> 0\\
 L^{ij}_{p} @>(\omega^{1}_{ij} \varepsilon+\omega^{2}_{ij} d\varepsilon) \eta>> L^{ij}_{0} \otimes \overline{\Omega}^{p}_{X \times k[\varepsilon]/k}(U_{ij})\otimes(\eta).
 \end{CD}
 \end{cases}
\]

 The following diagram (denoted $\beta_{ij}$)
\[
\begin{cases}
 \begin{CD}
   L^{ij}_{\bullet} @>>> O_{X}(U_{ij}) /(f_{i1}, \cdots, f_{ip}) @>>> 0\\
 L^{ij}_{p} @>g(\eta)(\omega^{1}_{ij} \varepsilon + \omega^{2}_{ij} d \varepsilon)>> L^{ij}_{0} \otimes \overline{\Omega}^{p}_{X \times k[\varepsilon]/k}(U_{ij}),
 \end{CD}
 \end{cases}
\]
defines an element in $Ext^{p}(O_{X}(U_{ij})/(f_{i1}, \cdots, f_{ip}), \overline{\Omega}^{p}_{X \times k[\varepsilon]/k}(U_{ij}))$, whose limit $[\beta_{ij}]$ is in $\Gamma(U_{ij}, \mathcal{H}_{Y}^{p}(\overline{\Omega}^{p}_{X \times k[\varepsilon]/k}))$. The 1-cocycle formed by $[\beta_{ij}]$ defines a class of $H^{1}(X, \mathcal{H}^{p}_{Y}(\overline{\Omega}^{p}_{X \times k[\varepsilon]/k})))$, which is $\delta_{e} \circ (\mathbb{T}^{(p)}(A),1)(Y',g)$.

 Since the diagram $\beta_{ij}$ agrees with the diagram $\overline{\beta}^{k[\varepsilon]}_{ij}$ (on page 29), it follows that
 \[
  \pi \circ v_{e}(Y', g)=\delta_{e} \circ (\mathbb{T}(A),1)(Y',g).
 \]

\end{proof}

In fact,  the map $\delta_{1}$ is trivial. By Lemma \ref{Lemma: HH-surjective}, for any element $\alpha_{A}=(\alpha^{A}_{i})_{i} \in \mathbb{HH}^{(p)}(A)$, one can lift $\alpha^{A}$ to $\alpha^{B} \in \mathbb{HH}^{(p)}(B)$. The restriction of $\alpha^{B}$ on each $U_{i}$, denoted $\alpha^{B}_{i}$, is a lifting of $\alpha^{A}_{i}$. Over the intersection $U_{ij}$, it is obvious that $\alpha^{B}_{i}-\alpha^{B}_{j}=0$. This shows that the map $\delta_{1}$ is trivial. Consequently, the map $\delta_{e}$ (4.6) is a trivial map. 

The commutative diagram of Lemma \ref{Lemma: Main Lemma-Commu-Diag} yields that the (set-theoretic) composition $\pi \circ v_{e}$
\[
\mathrm{\mathbb{H}ilb}(A)\times \check{(\eta)} \xrightarrow{v_{e}} H^{1}(Y, N_{Y/X}) \xrightarrow{\pi}   H^{1}(X, \mathcal{H}^{p}_{Y}(\overline{\Omega}^{p}_{X \times k[\varepsilon]/k}))
\]
is trivial.

By the description of the universal obstruction space $O_{\mathrm{\mathbb{H}ilb}}$ of the functor $\mathrm{\mathbb{H}ilb}$ (see Remark \ref{Remark: Universal Obs-Space}), $(e,Y')$ defines an element $[(e,Y')] \in O_{\mathrm{\mathbb{H}ilb}}$, where $e$ is the principal small extension (1.1) and $Y' \in \mathrm{\mathbb{H}ilb}(A)$.
The triviality of $\pi \circ v_{e}$ shows that the image of $[(e,Y')]$ under the composition of vector spaces
\[
O_{\mathrm{\mathbb{H}ilb}} \xrightarrow{[v_{e}]} H^{1}(Y, N_{Y/X}) \xrightarrow{\pi}   H^{1}(X, \mathcal{H}^{p}_{Y}(\overline{\Omega}^{p}_{X \times k[\varepsilon]/k}))
\]
is trivial, where $[v_{e}]$ is the linear monomorphism in Corollary \ref{Corollary: Obs HIlb is VS}. This immediately implies that 
\begin{theorem} \label{Theorem: Semi-reg obs space}
The image of $[(e,Y')]$ under the composition
\begin{align*}
&O_{\mathrm{\mathbb{H}ilb}} \xrightarrow{[v_{e}]} H^{1}(Y, N_{Y/X}) \xrightarrow{\pi} H^{1}(X, \mathcal{H}^{p}_{Y}(\overline{\Omega}^{p}_{X \times k[\varepsilon]/k})) \\
& \xrightarrow{\rfloor} H^{1}(X, \mathcal{H}^{p}_{Y}(\Omega^{p-1}_{X/k})) \xrightarrow{L} H^{p+1}(X, \Omega^{p-1}_{X/k})
\end{align*}
is trivial, where $L \circ  \rfloor \circ \pi$ is the semi-regularity map (see Lemma \ref{Lemma: Semi-regularity-Agree}).

\end{theorem}
 This gives a different proof of Theorem \ref{t:Bloch-semi-conj}.

\section{Appendix: Deformation functor}

Let $F$ be a functor of Artin rings. There is a natural map (see (2.6) on page 10)
\[
S: F(B \times_{A} C) \to F(B) \times_{F(A)} F(C).
\]

\begin{definition} [see Definition 3.4 of \cite{Manetti2005}]  \label{Definition: Deformation functor}
A functor of Artin rings $F$ is called a deformation functor if :
\begin{itemize}
\item [$\mathrm{(1)}$] Map $S$ is surjective if $C \to A$ is surjective,
\item [$\mathrm{(2)}$] Map $S$ is bijective if $A=k$.
\end{itemize}

\end{definition}

The local Hilbert functor $\mathrm{\mathbb{H}ilb}$ is an example of deformation functor. We explain that the functor of Artin rings $\mathbb{HH}^{(p)}$ in Definition \ref{Definition: Functor-HH-final} is not a deformation functor.

Let $A=k$, $B=k[x]/(x^{2})$ and $C=k[y]/(y^{2})$, then $B \times_{A} C= \{ (a+bx, a+cy)| a \in k, b \in k, c \in k \} $ and there is an isomorphism of $k$-algebras\footnote{The artinian $k$-algebra $k[x, y]/(x^{2}, xy, y^{2})$ had been studied by experts, for example, see Bloch \cite{Bloch3} (example 6.2.7), Manetti \cite{Manetti2005} (section 5) and Stienstra \cite{St} (the proof of Propsition 2.4).}
\[
B \times_{A} C \cong k[x, y]/(x^{2}, xy, y^{2}),
\] which implies that $\Omega^{1}_{B \times_{A} C /k}$ is a 3-dimensional $k$-vector space.

Since $\Omega^{1}_{A/k}=0$, $\Omega^{1}_{B/k} \times_{\Omega^{1}_{A/k}}\Omega^{1}_{C/k}=\Omega^{1}_{B/k} \times \Omega^{1}_{C/k}$ is a 2-dimensonal $k$-vector space. It is obvious that the map $\Omega^{1}_{B \times_{A} C /k} \to \Omega^{1}_{B/k} \times_{\Omega^{1}_{A/k}}\Omega^{1}_{C/k}$ is not bijective. Let $p=1$ in Lemma \ref{Lemma: HH-is localCoh}, one checks that the map 
\[
S: \mathbb{HH}^{(1)}(B \times_{A} C) \to \mathbb{HH}^{(1)}(B) \times_{\mathbb{HH}^{(p)}(A)}\mathbb{HH}^{(1)}(C)
\]
is not bijective. Hence, the functor $\mathbb{HH}^{(p)}$ does not satisfy the condition (2) in Definition \ref{Definition: Deformation functor}, it is not a deformation functor.

\textbf{Acknowledgments.} The author thanks Spencer Bloch \cite{Bloch2} for sharing his ideas and thanks him for comments on a preliminary version of this paper. He also thanks Jerome William Hoffman, Luc Illusie, Kefeng Liu and Chao Zhang for discussions, and thanks Shiu-Yuen Cheng and Bangming Deng for encouragement.

\end{document}